\documentclass{amsart} 

\usepackage{amssymb,latexsym,amsxtra,amscd,ifthen}
\usepackage{verbatim}
\usepackage[arrow,matrix,cmtip]{xy}

\usepackage{soul}
\usepackage[normalem]{ulem}

\providecommand{\qedhere}{}
\newboolean{ams2.0}
\ifthenelse{\equal{\qedhere}{}}%
	{\setboolean{ams2.0}{false}}%
	{\setboolean{ams2.0}{true}}%

\providecommand{\pdfoutput}{}
\ifthenelse{\equal{\pdfoutput}{}}%
			{}%
			{\usepackage[pdftex]{hyperref}}

\DeclareMathOperator{\GL}{GL}
\DeclareMathOperator{\Proj}{Proj}
\DeclareMathOperator{\proj}{proj}
\DeclareMathOperator{\Ext}{Ext}
\DeclareMathOperator{\Hom}{Hom}
\DeclareMathOperator{\calHom}{\mathcal{H}\mathit{om}}

\DeclareMathOperator{\Ker}{Ker}
\DeclareMathOperator{\Coker}{Coker}
\DeclareMathOperator{\Image}{Im}  
 
\DeclareMathOperator{\cd}{cd}

\DeclareMathOperator{\GKdim}{GKdim}
\DeclareMathOperator{\Supp}{Supp}
\DeclareMathOperator{\red}{red}
\DeclareMathOperator{\Pic}{Pic}
\DeclareMathOperator{\Spec}{Spec}
\DeclareMathOperator{\bSpec}{\mathbf{Spec}}

\DeclareMathOperator{\character}{char}

\DeclareMathOperator{\interior}{Int}

\DeclareMathOperator{\coh}{coh}
\DeclareMathOperator{\gr}{gr}

\DeclareMathOperator{\tors}{tors}

\DeclareMathOperator{\qgr}{qgr}

\newcommand{\ov}{\overline}

\newcommand{\op}{\text{op}}
\newcommand{\id}{\text{id}}

\newcommand{\scrm}{{\mathfrak m}}

\newcommand{\ZZ}{{\mathbb Z}}
\newcommand{\QQ}{{\mathbb Q}}
\newcommand{\RR}{{\mathbb R}}
\newcommand{\CC}{{\mathbb C}}
\newcommand{\NN}{{\mathbb N}}
\newcommand{\PP}{{\mathbb P}}
\newcommand{\LL}{{\mathcal L}}
\newcommand{\F}{{\mathcal F}}
\newcommand{\G}{{\mathcal G}}

\newcommand{\calC}{{\mathcal C}}

\newcommand{\calH}{{\mathcal H}}
\newcommand{\calI}{{\mathcal I}}
\newcommand{\calK}{{\mathcal K}}
\newcommand{\calN}{{\mathcal N}}
\newcommand{\calM}{{\mathcal M}}
\newcommand{\calP}{{\mathcal P}}
\newcommand{\OO}{{\mathcal O}}

\newcommand{\gs}{{\sigma}}

\numberwithin{equation}{section}

{\theoremstyle{plain}%
  \newtheorem{theorem}[equation]{Theorem}
  \newtheorem{corollary}[equation]{Corollary}
  \newtheorem{proposition}[equation]{Proposition}
  \newtheorem{lemma}[equation]{Lemma}%

{\theoremstyle{remark}

\newtheorem{remark}[equation]{Remark}
}
{\theoremstyle{definition}
\newtheorem{definition}[equation]{Definition}
\newtheorem{example}[equation]{Example}
\newtheorem{notation}[equation]{Notation}
}

\begin{document}

\title{Ample filters of invertible sheaves}

\author{Dennis S. Keeler}
     \thanks{ 
     Partially supported by an NSF Postdoctoral Research Fellowship.}
\address{Dept. of Mathematics \\ Miami University \\ Oxford, OH 45056 }
\email{keelerds@miamioh.edu}
\urladdr{http://www.users.miamioh.edu/keelerds/}
 

	{}%

%
%

\begin{abstract}  
Let $X$ be a scheme, proper over a commutative noetherian ring $A$.
We introduce the concept of an ample filter of invertible sheaves on $X$
and generalize the most important equivalent criteria
for ampleness of an invertible sheaf.
We also prove the Theorem of the Base for $X$
and generalize Serre's Vanishing Theorem.
We then generalize results for twisted homogeneous coordinate
rings which were previously known only when $X$ was projective
over an algebraically closed field. Specifically,
we show that the concepts of left and right $\sigma$-ampleness are
equivalent and that the associated twisted homogeneous coordinate ring
must be noetherian.
\end{abstract}

\maketitle


\section{Introduction} \label{S:intro}

Ample invertible sheaves are central to projective algebraic geometry.
Let $A$ be a commutative noetherian ring and let
$R$ be a commutative $\NN$-graded $A$-algebra, finitely generated in degree one.
Then ample invertible sheaves allow a geometric description
of $R$, by expressing $R$ as a homogeneous coordinate ring (in sufficiently high degree).
Further,  via the Serre Correspondence Theorem,
there is an equivalence between the category of coherent sheaves
on the scheme $\Proj R$ and the category of tails of graded $R$-modules
\cite[Exercise~II.5.9]{Ha1}. 

To prove that Artin-Schelter regular algebras of dimension $3$ (generated
in degree one) are noetherian, \cite{ATV} studied \emph{twisted homogeneous
coordinate rings} of elliptic curves (over a field).
In \cite{AV}, a more thorough and general examination of twisted homogenous coordinate
rings was undertaken, replacing the elliptic curves of \cite{ATV} with
any commutative projective scheme over a field.
Such a  ring $R$ is called twisted because it depends not only on
a  projective scheme $X$ and an invertible sheaf $\LL$, but also on an
automorphism $\gs$ of $X$, which causes the multiplication in $R$
to be noncommutative. We sometimes denote such a ring as $B(X, \gs, \LL)$.
When $\LL$ is \emph{right $\gs$-ample}, the ring $R$ is right noetherian and
 the category of tails of right $R$-modules is still equivalent 
to the category of coherent sheaves
on $X$.
When $\gs$ is the identity automorphism, the commutative theory is
recovered and $\gs$-ampleness reduces to the usual ampleness.
We will review $\gs$-ample invertible sheaves
and twisted homogeneous coordinate rings in $\S$\ref{S:rings}
of this paper.

A priori, there are separate definitions for right $\gs$-ampleness
and left $\gs$-ampleness.
However, when working with projective schemes
over an algebraically closed field, \cite{Ke1} shows
that right and left $\gs$-ampleness are equivalent.
In $\S$\ref{S:rings}, we generalize the results of 
\cite{AV} and \cite{Ke1} 
 to
the case of a scheme $X$ which is proper over a commutative noetherian ring $A$.
After summarizing definitions and previously known results,
 we prove
\begin{theorem}
Let $A$ be a commutative noetherian ring,
let $X$ be a proper scheme over  $A$, let $\gs$ be an automorphism  of $X$, and
let $\LL$ be an  invertible sheaf on $X$.
Then $\LL$ is right $\gs$-ample if and only if $\LL$
is left $\gs$-ample. If $\LL$ is $\gs$-ample,
then the twisted homogeneous coordinate ring $B(X, \gs, \LL)$
is noetherian.
\end{theorem}

To prove some of the results of \cite{AV} and \cite{Ke1}, the concept of $\gs$-ampleness
was not general enough. One can define the ampleness of a sequence of invertible
sheaves \cite[Definition~3.1]{AV}. To study twisted \emph{multi-homogenous} coordinate
rings, \cite{Chan} considers the ampleness of a set of invertible sheaves
indexed by $\NN^n$. To achieve the greatest possible generality, we index
invertible sheaves by any filter.

A filter $\calP$ is a   partially ordered set  such that:
\[
\text{for all } \alpha, \beta \in \calP, \text{ there exists } \gamma \in \calP
\text{ such that } \alpha < \gamma \text{ and } \beta < \gamma.
\]
Let $X$ be a scheme, proper over $\Spec A$, where $A$ is a (commutative) noetherian ring.
If a set of invertible sheaves is indexed by a filter, 
then we will call that
set a \emph{filter of invertible sheaves}. An element of such a 
filter will be denoted $\LL_\alpha$ for $\alpha \in \calP$. The indexing filter
$\calP$ will usually not be named.

\begin{definition}
Let $A$ be a commutative noetherian ring, and
let $X$ be a proper scheme over  $A$. Let $\calP$ be a filter.
A filter of invertible sheaves $\{ \LL_\alpha \}$ on $X$, with $\alpha \in \calP$,
 will be called an 
\emph{ample filter} if for all coherent sheaves $\F$ on $X$, there exists
$\alpha_0$ such that
\[
H^q(X, \F \otimes \LL_\alpha) = 0, \quad q > 0, \alpha \geq \alpha_0.
\]
If $\calP \cong \NN$ as filters, then an ample filter $\{ \LL_\alpha \}$
is called an \emph{ample sequence}.
\end{definition}

Of course if $\LL$ is an invertible sheaf, then
$\LL < \LL^{ 2} < \dots $ is an ample sequence if and only if
$\LL$ is an \emph{ample invertible sheaf}. 
It is well-known that the following conditions are equivalent in the
case of an ample invertible sheaf \cite[Proposition~III.5.3]{Ha1}. Our main result is 

\begin{theorem}\label{th:maintheo}
Let $X$ be a scheme, proper over a commutative noetherian ring $A$. Let $\{ \LL_\alpha \}$
be a filter of invertible sheaves. Then the following are
equivalent conditions on
$\{ \LL_\alpha \}$:
\begin{enumerate}
	\renewcommand{\theenumi}{A\arabic{enumi}}
	\item\label{th:main1} The filter $\{ \LL_\alpha \}$ is an ample filter.
	\item\label{th:main2} For all coherent sheaves $\F, \G$ with epimorphism 
		$\F \twoheadrightarrow \G$, there exists $\alpha_0$ such that the natural map
		\[
			H^0(X, \F \otimes \LL_\alpha ) \to H^0(X, \G \otimes \LL_\alpha)
		\]
		is an epimorphism for $\alpha \geq \alpha_0$.
	\item\label{th:main3} For all coherent sheaves $\F$, there exists $\alpha_0$ such that
		$\F \otimes \LL_\alpha$ is generated by global sections for $\alpha \geq \alpha_0$.
	\item\label{th:main4} For all invertible sheaves $\calH$, there exists $\alpha_0$ such that
	$\calH \otimes \LL_\alpha$ is an ample invertible sheaf for $\alpha \geq \alpha_0$.
\end{enumerate}
\end{theorem}
Note that there is no assumed relationship between the various $\LL_\alpha$, other
than their being indexed by a filter.

From condition \eqref{th:main4}, we see that if $X$ has an ample filter, then $X$ is
a projective scheme (see Corollary ~\ref{cor:projective}). A proper scheme $Y$
is divisorial if $Y$ has a so-called ample family of invertible
sheaves \cite[Definitions~ 2.2.4--5]{Illusie-divisorial}.
There exist divisorial schemes which are not projective; hence an ample family
does not have an associated ample filter in general. 
See Remark~ \ref{remark:divisorial} for some known
vanishing theorems on divisorial schemes.

Theorem~\ref{th:maintheo} is proven in $\S$\ref{S:ample-posets}. 
We must
first review previously known results involving ampleness and intersection
theory in $\S$\ref{S:review}. In particular, we review the general definitions
of numerical effectiveness and numerical triviality.

We then proceed to prove the Theorem of the Base in $\S$\ref{S:base}. 
 This was proven by Kleiman when $A$ was finitely
generated over a field \cite[p.~334, Proposition~3]{K}.
\begin{theorem}\label{th:base1}
{\upshape (See Thm.~\ref{th:base}) }
Let $X$ be a scheme, proper over a commutative noetherian ring $A$. Then
 $\Pic X$ modulo numerical equivalence is a finitely generated
free abelian group.
\end{theorem}
If two invertible sheaves $\LL, \LL'$ are numerically equivalent,
then $\LL$ is ample if and only if $\LL'$ is ample. Thus, one
may study ampleness via finite dimensional linear algebra.

After preliminary results on ample filters
in $\S$\ref{S:ample1},  we prove Serre's Vanishing Theorem
in new generality in $\S$\ref{S:serre}. The following was proven
by Fujita when
$A$ was an algebraically closed field \cite[$\S$5]{Fuj2}.

\begin{theorem}\label{th:general-serre1}
Let $A$ be a commutative noetherian ring,
let $X$ be a projective scheme over $A$, and let $\LL$ be an ample
invertible sheaf on $X$. For all coherent sheaves $\F$, there exists $m_0$ such that
\[
H^q(X, \F \otimes \LL^{ m} \otimes \calN) = 0
\]
for $m \geq m_0$, $q > 0$, and all  numerically effective invertible
sheaves $\calN$.
\end{theorem}

Finally, we note that before $\S$\ref{S:rings}, all rings in this paper are commutative.

\section{Ampleness and Intersection Theory}\label{S:review}
In this section, we review various facts about ampleness on a scheme which is
proper over a commutative noetherian ring $A$ or, more generally, proper over
a noetherian scheme $S$.
Most results are ``well-known,'' but proofs are not always easy to find
in the literature.
 We assume the
reader is familiar with various characterizations of ampleness which
appear in \cite{Ha1}. We mostly deal with 
noetherian schemes and proper morphisms, 
though some of the following propositions are true in
more generality than stated. 
All subschemes are assumed to be closed unless specified otherwise.

\begin{proposition} {\upshape \cite[II, 4.4.5]{EGA} }
Let $A$  be a  noetherian ring,
let $\pi\colon X \to \Spec A$ be a proper morphism, and
let 
$\LL$ be an invertible sheaf on $X$. 
If $U$ is an affine open subscheme of $\Spec A$ and $\LL$ is ample on $X$, 
then $\LL\vert_{\pi^{-1}(U)}$ is ample on $\pi^{-1}(U)$.
Conversely, if $\{ U_\alpha \}$ is  an affine open cover of $\Spec A$ and each
$\LL\vert_{\pi^{-1}(U_\alpha)}$ is ample on $\pi^{-1}(U_\alpha)$, then
$\LL$ is ample on $X$.\qed
\end{proposition}

This leads to the following definition of relative ampleness. 
\begin{definition} {\upshape \cite[II, 4.6.1]{EGA} }
Let $S$ be a noetherian scheme,
let $\pi\colon X \to S$ be a proper morphism,  and
let $\LL$ be an invertible sheaf on $X$. The sheaf $\LL$ is \emph{relatively ample}
for $\pi$ or \emph{$\pi$-ample} if there exists an affine open cover $\{ U_\alpha \}$
of $S$ such that $\LL\vert_{\pi^{-1}(U_\alpha)}$ is ample for each $\alpha$.
\end{definition}

\begin{proposition}\label{prop:ample-open}
Let $S$ be a noetherian scheme,
let $\pi\colon X \to S$ be a proper morphism,  and
let $\LL$ be an invertible sheaf on $X$.
If $U$ is an  open subscheme of $S$ and $\LL$ is $\pi$-ample, 
then $\LL\vert_{\pi^{-1}(U)}$ is  $\pi\vert_{\pi^{-1}(U)}$-ample.
Conversely, if $\{ U_\alpha \}$ is  an open cover of $S$ and each
$\LL\vert_{\pi^{-1}(U_\alpha)}$ is  $\pi\vert_{\pi^{-1}(U_\alpha)}$-ample, then
$\LL$ is  $\pi$-ample.\qed
\end{proposition}

\begin{remark}
If $\LL$ is $\pi$-ample, then some power $\LL^{ n}$ is 
relatively very ample for $\pi$ in the sense of Grothendieck
\cite[II, 4.4.2, 4.6.11]{EGA}.
Further a $\pi$-ample $\LL$ exists if and only if
$\pi$ is a projective morphism
\cite[II, 5.5.2, 5.5.3]{EGA}.
If $S$ itself has an ample invertible sheaf (for example when
$S$ is quasi-projective over an affine base), then 
these definitions of very ample invertible sheaf and projective morphism
are the same as those in \cite[p.~103, 120]{Ha1}, as shown in
\cite[II, 5.5.4(ii)]{EGA}. The definition of Grothendieck is useful because
 in many proofs
it allows an easy reduction to the case of  affine $S$.
\end{remark}

To find the connection between ampleness over general base rings and intersection
theory, we must examine how ampleness behaves on the fibers of a proper morphism.
We use the standard abuse of notation $X \times_{\Spec A} \Spec B = X \times_A B$.

\begin{proposition}\label{prop:fibers-ample} {\upshape \cite[$\mathrm{{III}_1}$, 4.7.1]{EGA} }
Let $S$ be a noetherian scheme,
let $\pi\colon X \to S$ be a proper morphism,  and
let $\LL$ be an invertible sheaf on $X$.
 Let $s \in S$ be a point
and $p\colon X_s = X \times_S  k(s) \to X$ be the natural projection
from the fiber. If $p^*\LL$ is ample on $X_s$, then there exists
an open neighborhood $U$ of $s \in S$ such that
$\LL\vert_{\pi^{-1}(U)}$ is  $\pi\vert_{\pi^{-1}(U)}$-ample.\qed
\end{proposition}

\begin{lemma}\label{lem:ZariskiSpace}
Let $S$ be a noetherian scheme
and let $U$ be an open subscheme which contains all closed points of $S$.
Then $U = S$.
\end{lemma}
\begin{proof}
Since $S$ is noetherian, it is a Zariski space. Thus, the closure
of any point of $S$ must contain a closed point \cite[Exercise~II.3.17e]{Ha1}.
So $S \setminus U$ must be empty.
\end{proof}

\begin{proposition}\label{prop:fibers-ample2}
Let $S$ be a noetherian scheme,
let $\pi\colon X \to S$ be a proper morphism,
and let $\LL$ be an invertible sheaf on $X$. Then $\LL$ is 
$\pi$-ample if and only if for every closed point $s \in S$
and natural projection $p_s\colon X_s \to X$, 
the invertible sheaf $p_s^* \LL$ is ample on $X_s$.
\end{proposition}
\begin{proof}
Suppose that $\LL$ is $\pi$-ample. We choose a closed point $s$
and an open affine neighborhood $s \in U \subset S$. Then $\LL\vert_{\pi^{-1}(U)}$
is ample on $\pi^{-1}(U) = X \times_S U$ by Proposition~\ref{prop:ample-open}.
Now since $s$ is closed, the fiber $X \times_S k(s)$ is a closed subscheme 
of $X \times_S U$.
So $p_s^* \LL$ is ample
\cite[Exercise~III.5.7]{Ha1}, as desired.

Conversely, suppose that $p_s^* \LL$ is ample for each closed point $s \in S$. 
Then by Proposition~\ref{prop:fibers-ample} and Lemma~\ref{lem:ZariskiSpace}, 
there is an open cover
$\{ U_\alpha \}$ of $S$ such that each $\LL\vert_{\pi^{-1}(U_\alpha)}$
is $\pi\vert_{\pi^{-1}(U_\alpha)}$-ample. So $\LL$
is  $\pi$-ample  by Proposition~\ref{prop:ample-open}.
\end{proof}

We now move to a discussion of 
 the intersection theory outlined in \cite[Chapter~VI.2]{KolRat}.
 We will soon see that the above proposition allows us to form a connection
 between this intersection theory and relative ampleness over a noetherian
 base scheme.

Let $S$ be a noetherian scheme, and let $\pi\colon X \to S$ be a proper morphism.
Let $\F$ be a coherent sheaf on $X$ with $\Supp \F$ proper over a $0$-dimensional
subscheme of $S$ and with $\dim \Supp \F = r$.
This theory defines intersection numbers $(\LL_1.\dots.\LL_n.\F)$ 
for the intersection of invertible sheaves $\LL_i$
on $X$ with  $\F$, where $n \geq r$. 
In the case $\F = \OO_Y$ for a closed subscheme $Y \subset X$, we also write
the intersection number as $(\LL_1.\dots.\LL_n.Y)$.
If all $\LL_i = \LL$, then we write $(\LL^{\bullet n}.Y)$. To avoid confusion
between self-intersection numbers and tensor powers, 
we will always write $\LL \otimes \dots \otimes \LL$
as $\LL^{\otimes n}$ when it appears in an intersection number.
(We use invertible sheaves instead of Cartier divisors since invertible
sheaves are more general \cite[Remark~II.6.14.1]{Ha1}.)

Note that
$\pi(\Supp \F)$ must be a closed subset of $S$ since $\pi$ is proper. 
According to \cite[Corollary~VI.2.3]{KolRat}, the intersection number
$(\LL_1.\dots.\LL_n.\F)$ may be calculated from
$(\LL_1.\dots.\LL_n.Y_i)$, where the $Y_i$ are the irreducible components
of the reduced scheme $\red \Supp \F$. Since $Y_i$ is irreducible and maps to
a $0$-dimensional subscheme of $S$,  we must have $\pi(Y_i)$ equal to
a closed point $s$ of $S$. Since $Y_i$ is reduced, there is a
unique induced proper morphism $f_{\red}\colon Y_i \to \Spec k(s)$, where
$k(s)$ is the residue field of $s$. Thus we may better understand this more 
general intersection theory by studying intersection theory over a field. 

The intersection theory of the case $S = \Spec k$, with $k$  algebraically closed, 
is the topic of the seminal paper \cite{K}. Many of the theorems of that paper are 
still valid in the case of $S$ equal to the spectrum of an arbitrary field $k$; the proofs can
either be copied outright or one can pass to $X \times_k \ov{k}$, where $\ov{k}$
is an algebraic closure of $k$. The most important proofs in this more general case
are included in \cite[Chapter~VI.2]{KolRat}. We will also need the following definitions
and lemmas.

\begin{definition}\label{def:pi-contracted}
Let $S$ be a noetherian scheme, and
let $\pi\colon X \to S$ be a proper morphism.
A closed subscheme $V \subset X$ is \emph{$\pi$-contracted} if
$\pi(V)$ is a $0$-dimensional (closed) subscheme of $S$.
If $S = \Spec A$ for a noetherian ring $A$, then we simply say that such a $V$
is  \emph{contracted}.
(This absolute notation will be justified by Proposition~\ref{prop:affine-base-not-relative}.)
\end{definition}

As stated above, if $V$ is irreducible and $\pi$-contracted, then $\pi(V)$
is a closed point. Recall that an integral curve is 
defined to be an integral scheme of dimension $1$ which is proper over a field.

\begin{definition}\label{def:nef}
Let $S$ be a noetherian scheme, and
let $\pi\colon X \to S$ be a proper morphism.
 An invertible sheaf $\LL$ on $X$ is
\emph{relatively numerically effective} for $\pi$ or \emph{$\pi$-nef}
(resp. \emph{relatively numerically trivial} for $\pi$ or \emph{$\pi$-trivial}) if
$(\LL.C) \geq 0$ (resp. $(\LL.C) = 0$) for all $\pi$-contracted integral curves $C$. 
If $S = \Spec A$ for a noetherian ring  $A$, then we simply say that such an $\LL$ is
\emph{numerically effective} or \emph{nef} (resp. \emph{numerically trivial}).
(This absolute notation will be justified by Proposition~\ref{prop:affine-base-not-relative}.)
\end{definition}

Obviously, $\LL$ is $\pi$-trivial if and only if $\LL$ is $\pi$-nef and
minus $\pi$-nef (i.e., $\LL^{ -1}$ is $\pi$-nef). Thus, the following propositions
regarding $\pi$-nef invertible sheaves have immediate corollaries for
$\pi$-trivial invertible sheaves, which we will not explicitly prove.

We need to study nef invertible sheaves because of their close connection
to ample invertible sheaves, as evidenced by the following propositions.
We begin with the Nakai criterion for ampleness. This is well-known 
when $X$ is proper over a field \cite[Theorem~VI.2.18]{KolRat}. The general
case follows via Proposition~\ref{prop:fibers-ample2}.

\begin{proposition}\label{prop:Nakai} {\upshape \cite[Thm.~1.42]{KM} }
Let $S$ be a noetherian scheme, let $\pi\colon X \to S$ be a proper morphism, 
and let $\LL$ be an invertible sheaf on $X$. Then $\LL$ is $\pi$-ample
if and only if $(\LL^{\bullet \dim V}.V) > 0$ for every
closed $\pi$-contracted  
 integral 
subscheme $V \subset X$.\qed
\end{proposition}

We have a similar proposition for nef invertible sheaves, following
from the case of $X$ proper over a base field \cite[Theorem~VI.2.17]{KolRat}.

\begin{proposition}\label{prop:pseudoample} {\upshape \cite[Thm.~1.43]{KM} }
Let $S$ be a noetherian scheme, let $\pi\colon X \to S$ be a proper morphism, 
and let $\LL$ be an invertible sheaf on $X$. Then $\LL$ is $\pi$-nef
(resp. $\pi$-trivial)
if and only if $(\LL^{\bullet \dim V}.V) \geq 0$ (resp. $(\LL^{\bullet \dim V}.V) = 0$) 
for every closed $\pi$-contracted  integral 
subscheme $V \subset X$.\qed
\end{proposition}

In fact, a stronger claim is true, namely

\begin{lemma}\label{lem:strong-pseudo} {\upshape \cite[p.~320, Thm.~1]{K}}
Given the hypotheses of Proposition~\ref{prop:pseudoample},
 if $\LL$ and $\calN$ are $\pi$-nef, 
then $(\LL^{\bullet i}.\calN^{\bullet \dim V - i}.V) \geq 0$ for every
closed $\pi$-contracted
 integral subscheme  and $i= 1, \dots, \dim V$.\qed
\end{lemma}
From this we easily see

\begin{proposition}\label{prop:ample-nef-connection}
Let $X$ be projective over a field $k$, and let $\LL$ be an ample invertible sheaf on $X$.
An invertible sheaf $\calN$ is numerically effective if and only if $\LL \otimes \calN^{ n}$
is ample for all $n \geq 0$.\qed
\end{proposition}

\begin{proposition}\label{prop:rel-ample-nef-connection}
Let $S$ be a noetherian scheme, let $\pi\colon X \to S$ be a proper morphism, 
and let $\LL$ be an invertible sheaf on $X$. 
An invertible sheaf $\calN$ is $\pi$-nef
if and only if $\LL \otimes \calN^{ n}$
is $\pi$-ample, for all $n \geq 0$.
\end{proposition}
\begin{proof}
Using Propositions~\ref{prop:fibers-ample2} and \ref{prop:ample-nef-connection},
we have the desired result by examining ampleness and numerical effectiveness
on the fibers over closed points.
\end{proof}

Ampleness over an affine base is an absolute notion 
\cite[Remark~II.7.4.1]{Ha1}, and this is also so for 
numerical effectiveness.

\begin{proposition}\label{prop:affine-base-not-relative}
Let $A_1, A_2$ be noetherian rings, and let
$\pi_i\colon X \to \Spec A_i$ be proper morphisms,
 $i=1,2$. Let $V$
be a closed subscheme of $X$. Then $V$ is
$\pi_1$-contracted if and only if $V$ is $\pi_2$-contracted. Thus
we may simply refer to such a $V$ as  contracted.

Let $\LL$ be an invertible sheaf
on $X$. Then $\LL$ is $\pi_1$-nef (resp. $\pi_1$-trivial) if and only
if $\LL$ is $\pi_2$-nef (resp. $\pi_2$-trivial). Thus we may simply
refer to such an $\LL$ as  numerically effective (resp. numerically
trivial).
\end{proposition}
\begin{proof}
Let $R = H^0(X, \OO_X)$. Since $X$ is proper over $A_1$, the ring $R$
is noetherian.
There are natural Stein factorizations of $\pi_i$ \cite[$\mathrm{III}_1$, 4.3.1]{EGA}
\[
\xymatrix@1{ X \ar[r]^-{\pi_i'} & \Spec R \ar[r]^{g_i} & \Spec A_i.}
\]
The $\pi_i'$ are both equal to the canonical morphism $f\colon X \to \Spec R$.
The $g_i$ are finite morphisms and $f$ is proper.
So a closed subscheme $V$ is $\pi_i$-contracted if and only if
it is $f$-contracted. 

Thus to prove the claim about $\LL$, 
we may assume the rings $A_i$ are fields
and $X$ is an integral curve. Then $\LL$ is  $\pi_i$-nef
if and only if $\LL^{-1}$ 
is not ample, and this is an absolute notion.
\end{proof}

We now examine the behavior of numerical effectiveness under pull-backs.
First we will need the following lemma.

\begin{lemma}\label{lem:surj-curve} 
Let $S$ be a noetherian scheme. Let
\[
\xymatrix{ X' \ar[r]^f \ar[rd]_{\pi'} & X \ar[d]^{\pi} \\
			& S }
\]
be a commutative diagram with proper morphisms $\pi, \pi'$ and
(proper) surjective morphism $f$. 
Let $C \subset X$ be a $\pi$-contracted integral curve. Then there exists a
$\pi'$-contracted integral
curve $C' \subset X'$ such that $f(C') = C$.
\end{lemma}
\begin{proof}
Note that $f$ is necessarily proper by \cite[Corollary~II.4.8e]{Ha1}.
If $\pi(C) = s$, we may replace $X, X', S$ by their fibers over
$\Spec k(s)$. Thus we are working over a field and this case
is as in \cite[p.~303, Lemma~1]{K}.
\end{proof}

\begin{lemma}\label{lem:surjective-nef} {\upshape (See \cite[p.~303, Prop.~1]{K})}
Let $S$ be a noetherian scheme. Let
\[
\xymatrix{ X' \ar[r]^f \ar[rd]_{\pi'} & X \ar[d]^{\pi} \\
			& S }
\]
be a commutative diagram with proper morphisms $\pi, \pi'$ and
(proper) morphism $f$.
Let $\LL$ be an invertible sheaf on $X$.
\begin{enumerate}
\item If $\LL$ is $\pi$-nef (resp. $\pi$-trivial), 
then $f^* \LL$ is $\pi'$-nef (resp. $\pi'$-trivial).
\item If $f$ is surjective and $f^* \LL$ is $\pi'$-nef (resp. $\pi'$-trivial), 
then $\LL$ is 
$\pi$-nef (resp. $\pi'$-trivial).
\end{enumerate}
\end{lemma}
\begin{proof}
By the definition of numerical effectiveness, we need only
examine the behavior of $\LL$ on the fibers of $\pi, \pi'$. Thus
we may assume $S = \Spec k$ for a field $k$.
 
For the first statement, let $C' \subset X'$ be an integral curve and let $C = f(C')$,
giving $C$ the reduced induced structure.
We have the Projection Formula $(f^* \LL.C') = (\LL. f_* \OO_{C'})$ 
\cite[Proposition~VI.2.11]{KolRat}.
Now 
\[
(\LL. f_* \OO_{C'}) = (\dim_{\OO_c} (f_* \OO_{C'})_c )(\LL.C)
\]
 where $c$ is the generic point
of $C$  \cite[Proposition~VI.2.7]{KolRat}. If $C$ is a point, then $(\LL.C) = 0$
and if $C$ is an integral curve, $(\LL.C) \geq 0$. So $f^* \LL$ is numerically effective.

For the second statement, let $C \subset X$ be an integral curve. By the previous lemma,
there is an integral curve $C' \subset X'$ such that $f(C') = C$. 
Since $f\vert_{C'}$ is surjective, $\dim_{\OO_c} (f_* \OO_{C'})_c$ is positive, so
the argument of the previous paragraph is reversible. That is,
\[
(\LL.C) = (f^* \LL.C')/(\dim_{\OO_c} (f_* \OO_{C'})_c) \geq 0.\qed
\]
\end{proof}

The preceding lemma is often used to reduce to the case of $X$ being projective  
over a noetherian scheme $S$. Given a proper morphism $\pi\colon X \to S$,
there exists a scheme $X'$ and a morphism $f\colon X' \to X$ such that
$f$ is projective and surjective, $\pi \circ f$ is projective, 
there exists a dense open subscheme $U \subset X$ such that
$f\vert_{f^{-1}(U)}\colon f^{-1}(U) \to U$ is an isomorphism, and
$f^{-1}(U)$ is dense in $X'$ \cite[II, 5.6.1]{EGA}. 
This $f\colon X' \to X$
is called a Chow cover.

Since \cite{K} works over an algebraically closed base field, the next lemma
is useful for applying results of that paper to the case of a general field.
In essence, it says the properties we are studying are preserved under base change.

\begin{lemma}\label{lem:field-extension}
Let $g \colon S' \to S$ be a morphism of noetherian schemes, let
$\pi\colon X \to S$ be a proper morphism, and let $\LL$ be an invertible sheaf on $X$.
Let $f\colon X \times_S S' \to X$ and $\pi'\colon X\times_S S' \to S'$ 
be the natural morphisms.
\begin{enumerate}
\item\label{lem:field-extension1} If $\LL$ is $\pi$-ample (resp. $\pi$-nef, $\pi$-trivial),
then $f^* \LL$ is 
$\pi'$-ample (resp. $\pi'$-nef, $\pi'$-trivial).
\item\label{lem:field-extension2} If $g(S')$ contains all closed points of $S$ 
and $f^* \LL$ is $\pi'$-ample
(resp. $\pi'$-nef, $\pi'$-trivial),
then $\LL$ is 
$\pi$-ample (resp. $\pi$-nef, $\pi$-trivial).
\end{enumerate}
\end{lemma}
\begin{proof}
For the first claim regarding ampleness, we may assume by Proposition~\ref{prop:ample-open}
that $S = \Spec A, S' = \Spec A'$ for noetherian rings $A, A'$.

Suppose that $\LL$ is ample. Then $\LL^{ n}$ is very ample for some $n$. We have the commutative
diagram
\[
\xymatrix{ X \times_A A' \ar@{^{(}->}[r]^-{i_{A'}} \ar@{->}[d]_f & \PP^m_{A'} \ar@{->}[d]_{f'} \\
		   X \ar@{^{(}->}[r]^-{i} 					   & \PP^m_A
		 }
\]
with $i, i_{A'}$ closed immersions. Now 
\[
f^* \LL^{ n} = f^*(i^* \OO_{\PP^m_A}(1)) = i_{A'}^* ({f'}^* \OO_{\PP^m_A}(1))
 = i_{A'}^* \OO_{\PP^m_{A'}}(1).
\]
Thus $f^* \LL^{ n}$ is very ample and $f^* \LL$ is ample.

Now suppose that $g(S')$ contains all closed points of $S$ and  $f^* \LL$ is ample. 
We wish to show $\LL$ is ample.
Let $s \in S$ be a closed point and let $s' \in S'$ be a closed point with $g(s')=s$.
Applying Proposition~\ref{prop:fibers-ample2} to $\LL$ and the first claim
of this lemma to $f^*\LL$, we may replace $S$ (resp. $S'$, $X$) with 
$\Spec k(s)$ (resp. $\Spec k(s')$, $X \times_S k(s)$).
Let $\F$ be a coherent sheaf on $X$.
Now $\Spec k(s') \to \Spec k(s)$ is faithfully flat,
so we have \cite[Proposition~III.9.3]{Ha1}
\[
H^q(X, \F \otimes \LL^m) \otimes_{k(s)} k(s') \cong H^q(X \times_S S', f^*\F \otimes f^*\LL^m) = 0
\]
for $m \geq m_0$, $q > 0$. Thus $\LL$ is ample.

For the claims regarding numerical effectiveness, we may replace $X$ with a Chow
cover via Lemma~\ref{lem:surjective-nef} and thus assume that $X$ is projective over $\Spec A$.
The claims then follow from Proposition~\ref{prop:rel-ample-nef-connection} and the
results above.
\end{proof}

Now we may generalize Proposition~\ref{prop:ample-open} as follows
for locally closed subschemes, i.e., a closed subscheme of an open subscheme.

\begin{corollary}\label{cor:nef-open}
Let $S$ be a noetherian scheme, let $\pi\colon X \to S$ be a proper morphism, 
and let $\LL$ be an invertible sheaf on $X$.
If $S_0$ is a locally closed subscheme of $S$ and $\LL$ is $\pi$-ample 
(resp. $\pi$-nef, $\pi$-trivial), 
then $\LL\vert_{\pi^{-1}(S_0)}$ is  $\pi\vert_{\pi^{-1}(S_0)}$-ample 
(resp. $\pi\vert_{\pi^{-1}(S_0)}$-nef, $\pi\vert_{\pi^{-1}(S_0)}$-trivial).
Conversely, if $\{ S_i \}$ is  a finite cover of locally closed subschemes of $S$ and each
$\LL\vert_{\pi^{-1}(S_i)}$ is  $\pi\vert_{\pi^{-1}(S_i)}$-ample 
(resp. $\pi\vert_{\pi^{-1}(S_i)}$-nef, $\pi\vert_{\pi^{-1}(S_i)}$-trivial), then
$\LL$ is  $\pi$-ample (resp. $\pi$-nef, $\pi$-trivial).
\end{corollary}
\begin{proof}
For the first claim, take $S' = S_0$ in 
Lemma~\ref{lem:field-extension}\eqref{lem:field-extension1}. For the second claim,
take $S' = \coprod S_i$ (the disjoint union of the $S_i$)
 in Lemma~\ref{lem:field-extension}\eqref{lem:field-extension2}.\footnote{The published
 version incorrectly writes $S' = \prod S_i$.}
\end{proof}

\begin{lemma}\label{lem:nef-change-scheme} 
Let $S, S'$ be noetherian schemes. Consider the commutative diagram
\[
\xymatrix{ X' \ar[r]^f \ar[d]_{\pi'} & X \ar[d]^\pi \\
			S' \ar[r]^g & S }
\]
with  proper morphisms $\pi, \pi'$. Let $\LL$ be an invertible sheaf
on $X$. 
\begin{enumerate}
\item If $\LL$ is $\pi$-nef (resp. $\pi$-trivial), then $f^*\LL$ is $\pi'$-nef
(resp. $\pi'$-trivial).
\item Further, suppose that 
for every $\pi$-contracted  integral curve $C \subset X$,
there is a $\pi'$-contracted  integral curve $C' \subset X'$
 such that $f(C') = C$.
 (For instance, suppose that $g = \id_S$ and $f$ is proper and surjective.)  
 If $f^*\LL$ is $\pi'$-nef (resp. $\pi'$-trivial), 
then $\LL$ is $\pi$-nef (resp. $\pi$-trivial).
\end{enumerate}
\end{lemma}
\begin{proof}
For the first statement, let $C' \subset X'$ be an integral curve such that
$\pi'(C') = s'$ is a closed point. Let $s = g(s')$. By the definition
of fibered products,  we have the induced commutative
diagram
\[
\xymatrix{ C' \ar[r]^-{f'} \ar[d]_{\pi'} 
	& X \times_{S} k(s') \ar[r]^-{f''} \ar[d] & X \times_S k(s) \ar[d]^\pi \\
		\Spec k(s') \ar@{=}[r] & \Spec k(s') \ar[r]^-{g''} & \Spec k(s).
		}
\]
Abusing notation, we replace $\LL$ with its restriction on $X \times_S k(s)$
and we need only show $(f')^* (f'')^* \LL$ is nef.
So assume $\LL$ is nef. By Lemma~\ref{lem:field-extension}, $(f'')^*\LL$ is nef.
And finally by Lemma~\ref{lem:surjective-nef}, $(f')^* (f'')^* \LL$ is nef.

Now assume the extra hypothesis of the latter part of the lemma. Let $C \subset X$
be any integral curve with $\pi(C) = s$ a closed point. Let $C' \subset X'$
be such that $f(C') = C$ and $\pi'(C') = s'$ is a closed point.
Then we abuse notation again, replacing $\LL$ with its restriction on $C$
and we have the induced commutative diagram
\[
\xymatrix{ C' \ar[r]^-{f'} \ar[d]_{\pi'} 
	& C \times_{k(s)} k(s') \ar[r]^-{f''} \ar[d]^{\pi''} & C \ar[d]^\pi \\
		\Spec k(s') \ar@{=}[r] & \Spec k(s') \ar[r]^-{g''} & \Spec k(s),
		}
\]
with $f = f'' \circ f'$.
Note that since $C, C'$ are proper curves over a field, they are projective
curves. Since $\pi'' \circ f'$ is projective and $\pi''$ is separated,
we have that $f'$ is projective \cite[Exercise~II.4.8]{Ha1}.
Thus $f'(C')$ is closed in $C \times_{k(s)} k(s')$.

Suppose that $(f')^* (f'')^* \LL$ is nef and
that $\LL$ is not nef, so that $(\LL.C) < 0$. Then
$\LL$ is minus ample. By Lemma~\ref{lem:field-extension}, $(f'')^*\LL$ is minus ample
on $C \times_{k(s)} k(s')$ and hence is minus ample on $f'(C')$, the closed image
subscheme. Since $f = f'' \circ f'$ is surjective and $C'$ is an integral curve,
$f'(C')$ must also be an integral curve. Since $f'$ is projective, no 
closed points of $C'$ can map to the generic point of $f'(C')$. The preimage
of any closed point of $f'(C')$ must also be a finite set, since finite sets
are the only proper closed subsets of $C'$. Thus, $f'$ is quasi-finite
and projective, hence finite \cite[Exercise~III.11.2]{Ha1}.
Thus $(f')^* (f'')^* \LL$ is minus ample on $C'$ \cite[Exercise~III.5.7]{Ha1},
contradicting that $(f')^* (f'')^* \LL$ is nef. So $\LL$ must be nef,
as desired.
\end{proof}

\section{The Theorem of the Base}\label{S:base}

Let $S$ be a noetherian scheme and let $\pi\colon X \to S$ be a proper morphism.
Two invertible sheaves $\LL, \LL'$ on $X$ are said to be \emph{relatively numerically
equivalent} if $(\LL.C) = (\LL'.C)$ for all $\pi$-contracted integral curves $C \subset X$.
If $S$ is affine, we say that $\LL$ and $\LL'$ are \emph{numerically equivalent}. 
We denote the equivalence by $\LL \equiv \LL'$.
Of course  $\LL \equiv 0$ exactly when $\LL$ is  relatively numerically trivial.
Define $A^1(X/S) = \Pic X/\equiv$. When $S = \Spec A$, we denote $A^1(X/S)$ as
$A^1(X/A)$ or simply $A^1(X)$, as the group depends only on $X$
by Proposition~\ref{prop:affine-base-not-relative}. 
The rank of $A^1(X)$ is called the
\emph{Picard number} of $X$ and is denoted $\rho(X)$.
Since the intersection numbers are multilinear,
the abelian group $A^1(X/S)$ is torsion-free.

The $A^1$ groups of two schemes may be related through the following
lemmas, which follow immediately from Lemmas~\ref{lem:nef-change-scheme}
and \ref{lem:field-extension}. We first generalize \cite[p.~334, Proposition~1]{K}.

\begin{lemma}\label{lem:nef-inject}
Let $S, S'$ be noetherian schemes. Let
\[
\xymatrix{ X' \ar[r]^f \ar[d]_{\pi'} & X \ar[d]^\pi \\
			S' \ar[r]^g & S }
\]
be a commutative diagram
with  proper morphisms $\pi, \pi'$. Then the pair $(f/g)$ induces a homomorphism
\[
(f/g)^*\colon A^1(X/S) \to A^1(X'/S').
\]
If for every $\pi$-contracted  integral curve $C \subset X$, 
there exists a $\pi'$-contracted  integral curve $C' \subset X'$ such that $f(C') = C$
(for instance if $g = \id_S$ and $f$ is proper and surjective), 
then $(f/g)^*$ is injective.\qed
\end{lemma}

\begin{lemma}\label{lem:faithfully-flat-inject}
Let $g \colon S' \to S$ be a  morphism of noetherian schemes, and let
$\pi\colon X \to S$ be a proper morphism.
Suppose that $g(S')$ contains all closed points of $S$. Then the natural map
$A^1(X/S) \to A^1(X \times_S S'/S')$ is injective.\qed
\end{lemma}

The Theorem of the Base says that $A^1(X/S)$ is finitely generated.
When $S = \Spec k$ with $k = \ov{k}$, this is \cite[p.~ 305, Remark~ 3]{K}.
The case of  an arbitrary base field $k$ can be
reduced to the algebraically closed case by Lemma~\ref{lem:faithfully-flat-inject}.

In \cite[p.~334, Proposition~3]{K}, the Theorem of the Base was proven when $S$ is
of finite-type over an algebraically closed field. 
We will follow much of this proof to prove Theorem~\ref{th:base}.
However, some changes must be made since normalization was used.
If $S$ is of finite-type over a field, the normalization of $S$ is still
a noetherian scheme, but
 there exists a noetherian (affine) scheme $S$ such that its normalization 
is not noetherian \cite[p.~127]{Eis}.
We will evade this difficulty via the following lemma. (We will
not use the claim regarding smoothness here, but it will be needed in $\S$\ref{S:serre}.)

\begin{lemma}\label{lem:geometrically-integral} 
\footnote{This lemma is incorrect. See the corrected Lemma~\ref{E-lem:geometrically-integral}. }
\sout{ Let $A$ be a noetherian domain with field of fractions $K$,
  let $X$ be an integral scheme
with  a projective, surjective morphism $\pi\colon X \to \Spec A$, 
and let $d$ be the dimension of the generic fiber of $\pi$. 
There exists non-zero $g \in A$,
a scheme $X'$, and a projective, surjective morphism $f \colon X' \to X \times_A A_g$
such that the composite morphism $\pi'\colon X' \to \Spec A_g$ is flat, projective,
and surjective,
and each fiber of $\pi'$ is geometrically integral, with
generic fiber dimension $d$.
Further, if $K$ is perfect, then one can assume that the morphism $\pi'$ is smooth. }
\end{lemma}
\begin{proof}
{ Let $X_0 = X \times_A K$ be the generic fiber of $\pi$. 
Using alteration of singularities \cite{deJong-Alterations}, 
we may find 
a regular integral $K$-scheme $\tilde{X}_0$ with projective, surjective morphism
$\tilde{X}_0 \to X_0$. Since an alteration is a generically finite morphism,
$\dim \tilde{X}_0 = d$. }

{ By \cite[Lemma~I.4.11]{J-Bertini}, there exists a finite extension field $K' \supset K$ such that
every irreducible component of $(\tilde{X}_0 \times_K K')_{\red{}}$
is geometrically integral.
Since $\tilde{X}_0 \times_K K' \to \tilde{X}_0$ is surjective
and finite, we may choose an irreducible component
$X_0'$ of $(\tilde{X}_0 \times_K K')_{\red{}}$ such that the
composite morphism $f_0\colon X_0' \to X_0$ is surjective
and projective. Since $K' \supset K$ is a finite extension,
$\dim X_0' = d$. }

{ Let $\pi_0' \colon X_0' \to \Spec K$ be the composite morphism.
Then $\pi_0'$ is certainly flat and projective. Further, if $K$ is perfect, then
 $\tilde{X}_0 \times_K K'$ is smooth over $\Spec K'$
and hence is a regular (and necessarily reduced) scheme \cite[Proposition~VII.6.3]{AltKle}.
Since two irreducible components of $\tilde{X}_0 \times_K K'$
cannot intersect \cite[Remark~III.7.9.1]{Ha1},  $X_0'$ is
a connected component (and hence a regular open subscheme) of $\tilde{X}_0 \times_K K'$.
Since $X_0'$ is algebraic over $K$, the morphism $\pi_0'$ is smooth. }

{ So now $\pi_0'$ and $f_0$ have all the desired properties of $\pi'$ and $f$.
Since $\pi_0'$ and $f_0$ are of finite type, we can find $g_1 \in A$ 
and an algebraic $A_{g_1}$-scheme $X_1'$ with finite type morphisms
\begin{equation}\label{eq:nice-morphisms}
\xymatrix{
X_1' \ar[r]^-{f_1} \ar[rd]_{\pi_1'} & X \times_A A_{g_1} \ar[d]^{\pi \times \id_{A_{g_1}}} \\
 & \Spec A_{g_1} }
\end{equation}
such that $\pi_1' \times_{A_{g_1}} \id_K \cong \pi_0'$ and $f_1 \times_{A_{g_1}} \id_K \cong f_0$.
We may shrink the base 
by taking the fibered product of \eqref{eq:nice-morphisms} with $\Spec A_g$ for
some multiple $g$ of $g_1$; we set $\pi' = \pi_1' \times \id_{A_g}$, $f = f_1 \times \id_{A_g}$
, $X' = X_1' \times A_g$.
By the Theorem of Generic Flatness, we may assume $\pi'$ is flat 
\cite[$\mathrm{IV}_3$, 8.9.4]{EGA}.
If $\{ A_\alpha \}$ is the inductive system 
of one element localizations of $A$, then $K = \varinjlim A_\alpha$.
So we may shrink the base to assume the resulting morphisms $\pi'$ and $f$
are projective and surjective \cite[$\mathrm{IV}_3$, 8.10.5]{EGA}.
Finally, we may assume all fibers  of $\pi'$ are geometrically integral
and (if $K$ is perfect) smooth over $k(s)$  \cite[$\mathrm{IV}_3$, 12.2.4]{EGA}.
So if $K$ is perfect, then $\pi'$ is smooth since $\pi'$ is flat
and all fibers are smooth over $k(s)$ \cite[Theorem~VII.1.8]{AltKle}. }
\end{proof}

We may now prove Theorem~\ref{th:base1}, which we state in the generality
of an arbitrary noetherian base scheme. 
First, a remark regarding fibers of flat morphisms.

\begin{remark}\label{remark:algebraic-closure-fiber}
Let $\pi\colon X \to S$ be a flat, projective morphism
of noetherian schemes. Let $s \in S$, and let $X_{\bar{s}} = X \times_S \ov{k(s)}$.
Let $\LL$ be an invertible sheaf on $X$, let $\LL_s = \LL\vert_{X_s}$
and let $\LL_{\bar{s}} = \LL_s \otimes_{k(s)} \ov{k(s)}$.
Since $\pi$ is flat, the Euler characteristic
$\chi(\LL_s)$ is independent of $s$ \cite[Lecture~7.9, Corollary~3]{MuCurves}.
Since $\Spec \ov{k(s)} \to \Spec k(s)$ is flat, 
$\chi(\LL_s) = \chi(\LL_s \otimes_{k(s)} \ov{k(s)})$.
\end{remark}
 
\begin{theorem}[Theorem of the Base]\label{th:base}
Let $S$ be a noetherian scheme and
let $\pi\colon X \to S$ be a proper morphism.
The torsion-free abelian group $A^1(X/S)$ is finitely generated.
\end{theorem}
\begin{proof}
Let $\{ U_\alpha \}$ be a finite open affine cover of $S$.
Using Corollary~\ref{cor:nef-open}, there is an induced injective
homomorphism 
\[ A^1(X/S) \hookrightarrow \oplus_\alpha A^1(\pi^{-1}(U_\alpha)/U_\alpha). \]
Thus we may assume $S = \Spec A$ is affine. By Lemma~ \ref{lem:nef-inject},
we may replace $X$ with a Chow cover, so we may assume $X$ is projective
over $\Spec A$. Let $X_i$ be the reduced, irreducible components of $X$. 
By Lemma~\ref{lem:nef-inject}, there is a natural monomorphism
$A^1(X) \hookrightarrow \oplus A^1(X_i)$, so we may assume $X$ is an
integral scheme. Also $A^1(X/\Spec A) = A^1(X/(\pi(X))_{\red})$, so
we may assume $A$ is a domain and $\pi$ is surjective.

We now proceed by noetherian induction on $S = \Spec A$. If $S = \emptyset$, then
$X = \emptyset$ and the theorem is trivial. Now assume that if $Y$ is any scheme 
which is projective over a proper closed subscheme of $S$, then $A^1(Y)$ is
finitely generated. The homomorphism 
$A^1(X) \to A^1(X \times_A [A/(g)]) \oplus A^1(X \times_A A_g)$
is injective, so we may replace $\Spec A$ with any affine open subscheme.
Then by Lemma~\ref{lem:nef-inject}, we may replace $X$ with the $X'$
of Lemma~\ref{lem:geometrically-integral} and assume $\pi$ is flat
and all fibers are geometrically integral.

Let $\eta$ be the generic point of $\Spec A$.
We claim that $\phi\colon A^1(X) \to A^1(X_\eta)$ is a
monomorphism. To see this, let $\LL \in \Ker \phi$, let $\calH$
be an ample invertible sheaf on $X$, and let $r = \dim X_\eta$.
By the definition of numerical triviality, Definition~\ref{def:nef}, 
we need to show that $\LL$ is numerically trivial on every closed fiber $X_s$
over  $s \in \Spec A$.
Equivalently by Lemma~\ref{lem:field-extension}, 
we need to show $\LL$ is numerically trivial on each $X_{\bar{s}} = X \times_A \ov{k(s)}$.

By Remark~\ref{remark:algebraic-closure-fiber}, 
the intersection numbers $(\LL_{\bar{s}}^{\bullet i}. \calH_{\bar{s}}^{\bullet r-i})$
are independent of $s$.
If $r = 0$, then any invertible sheaf on $X_{\bar{s}}$ is numerically trivial.
If $r = 1$, then $\LL_{\bar{s}}$ is numerically trivial if and only if $(\LL_{\bar{s}}) = 0$.
Finally, if $r \geq 2$, then $\LL_{\bar{s}}$ is numerically trivial if and only if
\[
(\LL_{\bar{s}}. \calH_{\bar{s}}^{\bullet r-1})
= (\LL_{\bar{s}}^{\bullet 2}. \calH_{\bar{s}}^{\bullet r-2}) = 0
\]
by \cite[p.~306, Corollary~3]{K}. Since $\LL_{\bar{\eta}}$ is numerically trivial,
we have shown that $\LL_{\bar{s}}$ is numerically trivial, independent of $s$,
and we are done.
\end{proof}

One of the most important uses of the Theorem of the Base is that it allows
the use of cone theory in studying ampleness and numerical effectiveness.
Now that we know that $A^1(X/S)$ is finitely generated, we may
prove Kleiman's criterion for ampleness in greater generality
than previously known. Since the original proofs rely
mainly on the abstract theory of cones in $\RR^n$, we may reuse
the original proofs in \cite[p.~323--327]{K}.
 A \emph{cone} $K$ is a subset of a real finite dimensional vector space $V$
such that for all $a \in \RR_{> 0}$ we have $aK \subset K$ and $K + K \subset K$.
The cone is \emph{pointed} if $K \cap -K = \{0\}$. The interior of 
a closed pointed cone $K$, $\interior K$, in the 
Euclidean topology of $V$, is also a cone, possibly empty.
\begin{lemma}\label{lem:interior} {\upshape \cite[p.~1209]{V}}
Let $K$ be a closed pointed cone in $\RR^n$.
Then $v \in \interior K$ if and only if, for all
$u \in \RR^n$, there exists $m_0$ such that
$u + mv \in K$ for all $m \geq m_0$.\qed
\end{lemma}

Now set $V(X/S) = A^1(X/S) \otimes \RR$. From Propositions~\ref{prop:Nakai} and
\ref{prop:pseudoample}, one sees that if $\LL, \LL'$ are invertible sheaves on $X$
with $\LL \equiv \LL' \in V(X/S)$, then $\LL$ is $\pi$-ample (resp. $\pi$-nef, $\pi$-trivial)
if and only if $\LL'$ has the same property. So for our purposes, we may use the notation
$\LL$ to represent an element of $\Pic X$ or $V(X/S)$ without confusion.

One can easily show that the cone $K$ generated by $\pi$-nef invertible sheaves
is a closed pointed cone. The cone $K^\circ$ generated by $\pi$-ample invertible
sheaves is open in the Euclidean topology \cite[p.~325, Remark~6]{K} and 
$K^\circ \subset \interior K$ by Lemma~\ref{lem:interior}.
We would like to say that $K^\circ = \interior K$ for all proper schemes,
but this is not true for some degenerate non-projective cases.
For example, in \cite[Exercise~III.5.9]{Ha1}, $\Pic X = 0$, so
$K^\circ = \emptyset \subsetneq \interior K = K = 0$.
We need a suitable generalization of projectivity.

\begin{definition}
Let $\pi\colon X \to S$ be a proper morphism over a noetherian scheme $S$.
The scheme $X$ is \emph{relatively quasi-divisorial} for $\pi$ if for every
 $\pi$-contracted
integral subscheme $V$ (which is not a point), there
exists an invertible sheaf $\LL$ on $X$ and an effective non-zero Cartier
divisor $H$ on $V$ such that $\LL\vert_Y \cong \OO_Y(H)$. If
$S$ is affine, then $X$ is \emph{quasi-divisorial}.
(This absolute notation is justified by Proposition~\ref{prop:affine-base-not-relative}.)
\end{definition}

If $\pi$ is projective, then $X$ is relatively quasi-divisorial; just
take $\LL$ to be relatively very ample for $\pi$. If $X$ is a regular integral scheme
 (or more generally $\QQ$-factorial), 
 then $X$ is relatively quasi-divisorial \cite[Pf. of Theorem~VI.2.19]{KolRat}.
 See Remark~\ref{remark:divisorial} for the definition of a divisorial scheme.

\begin{theorem}[Kleiman's criterion for ampleness]\label{th:Kleiman} 
Let $\pi\colon X \to S$ be a proper morphism over a noetherian scheme $S$
with $X$ relatively quasi-divisorial for $\pi$.
An invertible sheaf $\LL$ on $X$ is $\pi$-ample if and only if $\LL \in \interior K$.
More generally, $K^\circ = \interior K$.
\end{theorem}
\begin{proof}
By abstract cone theory, the cone generated by
the lattice points $\interior K \cap A^1(X/S)$ is equal to $\interior K$ 
\cite[p.~325, Remark~5]{K}. So we need only show the first claim.
We already know that if $\LL$ is $\pi$-ample, then $\LL \in \interior K$.

Now suppose that $\LL \in \interior K$. By Proposition~\ref{prop:Nakai}, we need
to show that for every  $\pi$-contracted integral subscheme $V$ of $X$,
we have $(\LL^{\bullet \dim V}.V) > 0$. Since $X$ is relatively
quasi-divisorial, there exists $\calH \in \Pic X$ such that $\calH \cong \OO_V(H)$
with $H$ a non-zero effective Cartier divisor.
Note that if $\dim V = 1$, then $(\calH.V) = (H) = \dim H^0(H,\OO_H) > 0$
\cite[p.~296, Proposition~1]{K}.

For sufficiently large $n$, we have 
$\LL^{ n} \otimes \calH^{ -1} \in K$ by Lemma~\ref{lem:interior}.
And so, $(\LL^{{\bullet} \dim V - 1}.\LL^{\otimes n} \otimes \calH^{ -1}.V) \geq 0$
by Lemma~\ref{lem:strong-pseudo}.
Thus
\[
n(\LL^{\bullet \dim V}.V) \geq (\LL^{\bullet \dim V - 1}.\calH.V) = (\LL^{\bullet \dim V - 1}.H) > 0
\]
by induction. Thus $\LL$ is ample.
\end{proof}

\section{Ample filters I}\label{S:ample1}

In this section we will collect a few preliminary propositions regarding ample filters
which are well-known in the case of an ample invertible sheaf.
The main goal of $\S$\ref{S:serre} will be to prove that a certain filter of invertible
sheaves is an ample filter; these propositions will allow for useful reductions
in that proof.

\begin{proposition}\label{prop:filter-reductions}
Let $X, Y$ be proper over a noetherian ring $A$. Let $\{ \LL_\alpha \}$
be a filter of invertible sheaves on $X$.
\begin{enumerate}
\item If $\{ \LL_\alpha \}$ is an ample filter on $X$, then $\{ \LL_\alpha\vert_Y \}$
is an ample filter on $Y$ for all closed subschemes $Y \subset X$.
\item The filter $\{ \LL_\alpha \}$ is an ample filter on $X$ if and only if
 $\{ \LL_\alpha\vert_{X_{\red}} \}$ is an ample filter on $X_{\red}$.
\item Suppose that $X$ is reduced. The filter $\{ \LL_\alpha \}$ is an ample filter on $X$ if and only if
for each irreducible component $X_i$, the filter
 $\{ \LL_\alpha\vert_{X_i} \}$ is an ample filter on $X_i$.
\item\label{pull-back4} Let $f\colon Y \to X$ be a finite morphism. 
If $\{ \LL_\alpha \}$ is an ample filter on $X$,
then $\{ f^*\LL_\alpha \}$ is an ample filter on $Y$.
\item\label{pull-back5} Let $f\colon Y \to X$ be a finite surjective morphism.
If $\{ f^*\LL_\alpha \}$ is an ample filter on $Y$,
then $\{ \LL_\alpha \}$ is an ample filter on $X$.
\end{enumerate}
\end{proposition}
\begin{proof}
The proof of each item is as in \cite[Exercise~III.5.7]{Ha1}.
\end{proof}

For the rest of this section we will switch from the proper case to
the projective case. However, this is not a real limitation in
regards to ample filters since we will see that if $X$ has
an ample filter, then $X$ is projective by Corollary~\ref{cor:projective}.

\begin{proposition}\label{prop:vanishing-ample-powers-enough}
Let $X$ be projective over a noetherian ring $A$, let
$\calH$ be an ample invertible sheaf on $X$, and let
$\{ \LL_\alpha \}$ be a filter of invertible sheaves on $X$.
The filter $\{ \LL_\alpha \}$ is an ample filter if and only if
for each $n > 0$, there exists $\alpha_0$ such that
\[
H^q(X, \calH^{-n} \otimes \LL_\alpha) = 0
\]
for $q>0$ and $\alpha \geq \alpha_0$.
\end{proposition}
\begin{proof}
Necessity is obvious. So assume
the vanishing of cohomology of $\calH^{-n} \otimes \LL_\alpha$.
We claim that
for each $q > 0$ and coherent sheaf $\F$ on $X$, 
there exists $\alpha_q$ such that
$H^{q'}(X, \F \otimes \LL_\alpha) = 0$ for $\alpha \geq \alpha_{q}$ and
$q' \geq q$.
We have a short exact sequence
\begin{equation}\label{eq:cokernel-of-ample}
0 \to \calK \to \oplus_{i=1}^p \calH^{-n} \to \F \to 0
\end{equation}
for some $n > 0$ and some coherent sheaf $\calK$. Then
there exists $\alpha_0$ such that
$H^q(X, \F \otimes \LL_\alpha) \cong H^{q+1}(X, \calK \otimes \LL_\alpha)$
for all $q > 0$ and $\alpha \geq \alpha_0$.
For $q+1 > \cd(X)$, we have $H^{q+1}(X, \calK \otimes \LL_\alpha) = 0$,
where $\cd(X)$ is the cohomological dimension of $X$ \cite[Exercise~III.4.8]{Ha1}.
So by descending induction on $q$, the claim is proved and so is the proposition.
\end{proof}

As in the case of ample invertible sheaves, some propositions are best
stated in the relative situation. So we make the following definition.

\begin{definition}
Let $S$ be a noetherian scheme, and let $\pi\colon X \to S$ be a proper
morphism.
A filter of invertible sheaves $\{ \LL_\alpha \}$ on $X$, with $\alpha \in \calP$,
 will be called a 
\emph{$\pi$-ample filter} if for all coherent sheaves $\F$ on $X$, there exists
$\alpha_0$ such that
\[
R^q\pi_*(\F \otimes \LL_\alpha) = 0, \quad q > 0, \alpha \geq \alpha_0.
\]
If $\calP \cong \NN$ as filters, then a $\pi$-ample filter $\{ \LL_\alpha \}$
is called a \emph{$\pi$-ample sequence}.
\end{definition}

We now state a partial generalization of Corollary~\ref{cor:nef-open}.

\begin{proposition}\label{prop:ample-filter-locally-closed}
Let $S$ be a noetherian scheme, let $\pi\colon X \to S$ be a proper morphism, 
and let $\{ \LL_\alpha \}$ be a filter of invertible sheaves on $X$.
If $S_0$ is a locally closed subscheme of $S$ and $\{ \LL_\alpha \}$ is a $\pi$-ample 
filter, 
then the filter $\{ \LL_\alpha\vert_{\pi^{-1}(S_0)} \}$ is  a $\pi\vert_{\pi^{-1}(S_0)}$-ample
filter.
\end{proposition}
\begin{proof}
The scheme $S_0$ is a closed subscheme of an open subscheme $U$ of $S$.
Let $\F$ be a coherent sheaf on $X \times_S S_0$. Then
 $i\colon X \times_S S_0 \hookrightarrow X \times_S U$
is a closed immersion and
$\F' = i_* \F$ is a coherent sheaf
on $X \times_S U$.
Further, $v\colon X \times_S U \hookrightarrow X$ is an open immersion
and
 there exists a coherent sheaf $\F''$ on $X$ such that
$v^*\F'' = \F'$ \cite[Exercise~II.5.15]{Ha1}.
We have the desired vanishing of higher direct images
after tensoring with $i^* v^* \LL_\alpha = \LL_\alpha\vert_{\pi^{-1}(S_0)}$
 since $u \colon U \hookrightarrow S$
is flat \cite[Proposition~III.9.3]{Ha1}
and $i$ is affine \cite[Exercise~III.4.1]{Ha1}.
\end{proof}

The following lemma will be useful in $\S$\ref{S:serre} for noetherian induction
on $X$. Its proof has similarities to the proof of
Proposition~\ref{prop:fibers-ample}.

\begin{lemma}\label{lem:open-around-closed}
Let $S$ be a noetherian scheme, let $\pi\colon X \to S$ be a projective morphism,
and let $\{ \LL_\alpha \}$ be a filter of invertible sheaves on $X$.
Let $T$ be a closed subscheme of $S$ and suppose that
$\{ \LL_\alpha\vert_{\pi^{-1}(T)} \}$ is a $\pi\vert_{\pi^{-1}(T)}$-ample
filter on $\pi^{-1}(T) = X \times_S T$.
For each coherent sheaf $\F$ on $X$, there exists $\alpha_0$
and open subschemes $U_\alpha \supset T$, 
such that
\[
R^q \pi_*(\F \otimes \LL_\alpha)\vert_{U_\alpha} = 0
\]
for $q > 0$ and $\alpha \geq \alpha_0$.
\end{lemma}
\begin{proof}
If $T = \emptyset$, then we may take $U_\alpha = \emptyset$ for all
$\alpha$. So by noetherian induction, we may assume
that for all proper closed subschemes $V \subset T$
and coherent sheaves $\F$ on $X$, there exists $\alpha_0$
and $U_\alpha \supset V$ such that
\[
R^q \pi_*(\F \otimes \LL_\alpha)\vert_{U_\alpha} = 0
\]
for $q > 0$ and $\alpha \geq \alpha_0$.
Further, if $U$ is an open subscheme of $S$ and $W = (U \cap T)_{\red}$, then
$\{ \LL_\alpha\vert_{\pi^{-1}(W)} \}$ is 
a $\pi\vert_{\pi^{-1}(W)}$-ample filter by Proposition~\ref{prop:ample-filter-locally-closed}.
So we may replace $S$ with an open subscheme $U$ as long as $U \cap T \neq \emptyset$.
Thus, we may assume $S = \Spec A$ and $T = \Spec A/I$.

Let $B = \gr_I(A) = \oplus_{j=0}^\infty I^j/I^{j+1}$. 
Then $\Spec B$ is called the normal cone $C_T S$ to $T$ in $S$
\cite[Appendix~B.6.1]{Fu}.
Let $R$ be a ring of polynomials over $A/I$ which surjects
onto $B$.
Then $u \colon \Spec R \to \Spec A/I$ is flat.

Now consider the commutative diagram
\[
\xymatrix{ 
		X \times_A R \ar[r]^g \ar[d] 
		& X \times_A A/I \ar@{^{(}->}[r]^-{i} \ar[d]
		& X \ar[d]_\pi \\
		\Spec R \ar[r]^u & \Spec A/I \ar@{^{(}->}[r] & \Spec A. }
\]
Let $\calH$ be an ample invertible sheaf
on $X$, and let $h = i \circ g$. Then by Lemma~\ref{lem:field-extension}, 
$h^* \calH$ is an ample invertible sheaf. Given $n > 0$, 
there exists $\alpha_0$ such that
\[
H^q(X \times_A R, h^*(\calH^{-n} \otimes \LL_\alpha)) =
H^q(X \times_A A/I, i^*(\calH^{-n} \otimes \LL_\alpha)) \otimes_{A/I} R = 0
\]
for $\alpha \geq \alpha_0$, since $u$ is flat and $\{ i^*\LL_\alpha \}$ 
is an ample filter on $X \times_A A/I$.
Thus $\{ h^*\LL_\alpha \}$ is an ample filter on $X \times_A R$
by Proposition~\ref{prop:vanishing-ample-powers-enough}.

Let $\calI = \pi^{-1}(\tilde{I})\OO_X$ be the sheaf
of ideals of $X \times_A A/I$ in $X$. Then there is a canonical embedding
\cite[Appendix~B.6.1]{Fu}
\[
C_{X \times_A A/I} X = \bSpec \bigoplus_{j=0}^\infty \calI^j/\calI^{j+1}
\stackrel{i'}{\hookrightarrow} X \times_A B = X \times_A C_T S
\]
and a closed embedding $i''\colon X \times_A B \hookrightarrow X \times_A R$.
Let $f =  h \circ i'' \circ i'$.
Then by Proposition~\ref{prop:filter-reductions}, $\{ f^* \LL_\alpha \}$
is an ample filter on $C_{X \times_A A/I} X$.
We set $Y = C_{X \times_A A/I} X$.

Let $\F$ be a coherent sheaf on $X$, and let
$\F' = \oplus_{j=0}^\infty \calI^j\F/\calI^{j+1}\F$.
Then $\F'$ is a quasi-coherent $\OO_X$-module with an obvious $f_*\OO_Y$-module
structure. Since $f$ is affine, there exists a quasi-coherent $\OO_Y$-module
$\G$ such that $f_*\G = \F'$ \cite[Exercise~II.5.17]{Ha1}.
But since $\F'$ is a coherent $f_*\OO_Y$-module, the $\OO_Y$-module
$\G$ is coherent.
So there exists $\alpha_1$ such that
\begin{align*}
\bigoplus_{j=0}^\infty H^q(X, \calI^j\F/\calI^{j+1}\F \otimes \LL_\alpha)
&\cong H^q(X, \F' \otimes \LL_\alpha) \\
&\cong H^q(X, f_*(\G \otimes f^*\LL_\alpha)) \\
&\cong H^q(Y, \G \otimes f^*\LL_\alpha) = 0
\end{align*}
for $q > 0$ and $\alpha \geq \alpha_1$, because $f$ is affine \cite[Exercise~III.4.1]{Ha1}
and $\{ f^*\LL_\alpha \}$ is an ample filter.

Thus $H^q(X, \calI^j\F/\calI^{j+1}\F \otimes \LL_\alpha) = 0$
for all $q>0$, $j \geq 0$ and $\alpha \geq \alpha_1$.
Now consider the short exact sequence
\[
0 \to \F/\calI^j\F \to \F/\calI^{j+1}\F \to \calI^j\F/\calI^{j+1}\F \to 0
\] 
which yields an exact sequence
\[
H^q(X, \F/\calI^j\F \otimes \LL_\alpha) 
\to H^q(X, \F/\calI^{j+1}\F \otimes \LL_\alpha)
\to H^q(X, \calI^j\F/\calI^{j+1}\F \otimes \LL_\alpha).
\]
The last term is $0$ for all $j \geq 0$ and the first term is $0$ for
$j = 1$. Then by induction, the middle term is $0$ for all $q>0$, $j \geq 0$,
and $\alpha \geq \alpha_1$.

Let $\hat{A}$ be the $I$-adic completion of $A$.
By \cite[$\mathrm{III}_1$, 4.1.7]{EGA},
\[
 H^q(X, \F \otimes \LL_\alpha) \otimes_A \hat{A} \cong
\varprojlim_j H^q(X, \F/\calI^{j+1}\F \otimes \LL_\alpha  ) = 0
\]
for $q > 0$ and $\alpha \geq \alpha_1$. Thus for each
$\alpha \geq \alpha_1$, there exists $a_\alpha \in 1 + I$
such that $H^q(X, \F \otimes \LL_\alpha) \otimes_A A_{a_\alpha} = 0$ \cite[Theorem~10.17]{AM}. 
Since
$a_\alpha \in 1 + I$, the open subscheme $\Spec A_{a_\alpha}$
contains the closed subscheme $\Spec A/I$, as desired.
\end{proof}

In Proposition~\ref{lem:open-around-closed} we generalized
 the first half of Corollary~\ref{cor:nef-open}.
 We now generalize the second half,
which will allow the use of noetherian induction in $\S$\ref{S:serre}
to show that a certain filter is an ample filter.

\begin{proposition}\label{prop:ample-cover-implies-ample-filter}
Let $S$ be a noetherian scheme, let $\pi\colon X \to S$ be a proper morphism, 
and let $\{ \LL_\alpha \}$ be a filter of invertible sheaves on $X$.
 If $\{ S_i \}$ is  a finite cover of locally closed subschemes of $S$ and each
$\{ \LL_\alpha \vert_{\pi^{-1}(S_i)} \}$ is a $\pi\vert_{\pi^{-1}(S_i)}$-ample 
filter, then $\{ \LL_\alpha \}$ is a $\pi$-ample filter.
\end{proposition}
\begin{proof}
The schemes $S_i$ are closed subschemes of open subschemes $U_i$ of $S$.
Let $\F$ be a coherent sheaf on $X$.
Then by Lemma~\ref{lem:open-around-closed}, for each $i$, there exists $\alpha_i$
and open subschemes $U_{i,\alpha}$ of $U_i$ such that
$S_i \subset U_{i,\alpha}$ and \cite[Corollary~III.8.2]{Ha1}
\[
R^q (\pi\vert_{\pi^{-1}(U_{i})})_*
((\F \otimes \LL_\alpha)\vert_{\pi^{-1}(U_{i})})\vert_{U_{i,\alpha}}
= R^q \pi_*(\F \otimes \LL_\alpha)\vert_{U_{i,\alpha}}
= 0
\]
for $q > 0$ and $\alpha \geq \alpha_i$. Let $\alpha_0 \geq \alpha_i$ for all $i$.
Then for $\alpha \geq \alpha_0$, the open subschemes $U_{i, \alpha}$ cover $S$,
so $R^q\pi_*(\F \otimes \LL_\alpha) = 0$, as desired.
\end{proof} 

We will also need a theorem for global generation of $\F \otimes \LL_\alpha$.
We do so through the concept of $m$-regularity, which is most often studied when
$X$ is projective over a field. However, most of the proofs are still valid
when $X$ is projective over a noetherian ring $A$.

\begin{definition}
Let $X$ be a projective scheme over a noetherian ring $A$, and let
 $\OO_X(1)$ be a very ample invertible
sheaf. A coherent sheaf $\F$ is said to be \emph{$m$-regular}
(with respect to $\OO_X(1)$)
if $H^q(X, \F(m-q) ) = 0$ for $q > 0$.
\end{definition}

\begin{proposition}\label{prop:regularity}
Let $X$ be a projective scheme over a noetherian ring $A$, and let
 $\OO_X(1)$ be a very ample invertible
sheaf. If a coherent sheaf $\F$ on $X$ is $m$-regular, then for $n \geq m$
\begin{enumerate}
\item\label{regular1} $\F$ is $n$-regular,
\item\label{regular2} The natural map 
$H^0(X, \F(n) ) \otimes_A H^0(X, \OO_X(1) ) \to H^0(X, \F(n+1) )$
is surjective, and
\item\label{regular3} $\F(n)$ is generated by global sections.
\end{enumerate}
\end{proposition}
\begin{proof}
Let $j\colon X \hookrightarrow \PP^r_A$ be the closed immersion defined by $\OO_X(1)$.
Then for all $n$ and $q \geq 0$, 
$H^q(X, \F \otimes \OO_X(n) ) \cong H^q(\PP^r, (j_*\F) \otimes \OO_{\PP^r}(n) )$.
Thus $\F$ is $n$-regular if and only if $j_* \F$ is $n$-regular.
So \eqref{regular1} holds on $X$ if it holds on $\PP^r$. Further, we have
the commutative diagram
\[
\xymatrix
{ 
{ H^0( \PP^r, j_* \F(n) ) \otimes H^0( \PP^r, \OO_{\PP^r}(1)) } \ar[r] \ar[d]  
	&  H^0(\PP^r, j_* \F(n+1))  \\ 
	 {H^0( \PP^r, j_* \F(n) ) \otimes H^0( \PP^r, j_*\OO_X(1)), \ar[ur]} 
	}
\]
so if \eqref{regular2} holds on $\PP^r$, it holds on $X$.
So for \eqref{regular1} and \eqref{regular2} 
we have reduced to the case $X = \PP^r$ and this is \cite[Theorem~2]{Oo}.

The proof of \eqref{regular3} 
proceeds as in \cite[p.~307, Proposition~1(iii)]{K}, keeping in
mind the more general situation. Let $f\colon X \to \Spec A$ be the structure morphism.
A coherent sheaf $\G$ is generated by global sections if and only if
the natural morphism $f^* f_* \G = f^* \widetilde{H^0(X, \G)} \to \G$ is surjective
\cite[Theorem~III.8.8]{Ha1}.
We have a commutative diagram
\[
\xymatrix{ f^*\widetilde{H^0(\F(n))} \otimes f^*\widetilde{H^0(\OO_X(1))}
			\ar[r]^-{\alpha_n} \ar[d] & f^*\widetilde{H^0(\F(n+1))} \ar[d]_{\beta_{n+1}} \\
			f^*\widetilde{H^0(\F(n))} \otimes \OO_X(1) \ar[r]^-{\beta_n \otimes 1} & \F(n+1).
		}
\]
By \eqref{regular2}, $\alpha_n$ is surjective for $n \geq m$. Also, there
exists $n_1 \geq m$  such that
$\F(n+1)$ is generated by global section for $n\geq n_1$, and so $\beta_{n+1}$ is also
surjective for $n \geq n_1$. This implies that $\beta_n \otimes 1$ (and hence $\beta_n$) 
is surjective
for $n \geq n_1$. Descending induction on $n$ gives that $\beta_n$ is surjective for $n \geq m$, 
as desired.
\end{proof}

\begin{corollary}\label{cor:general-serre-global-sections}
Let $X$ be a projective scheme over a commutative noetherian ring $A$,  let $\OO_X(1)$
be a very ample
invertible sheaf on $X$, let $\{ \LL_\alpha \}$ be an ample filter
on $X$, and let $\F$ be a coherent sheaf on $X$.
There exists $\alpha_0$ such that
\begin{enumerate}
\item The natural map 
$H^0(X, \F \otimes \LL_\alpha) \otimes H^0(X, \OO_X(1))
\to H^0(X, \F \otimes \LL_\alpha \otimes \OO_X(1))$ is surjective and
\item $\F \otimes \LL_\alpha$ is generated by global sections
\end{enumerate}
for $\alpha \geq \alpha_0$.
\end{corollary}
\begin{proof}
Find $\alpha_0$ such that 
$H^q(X, \F \otimes \LL_\alpha \otimes \OO_X(-q)) = 0$ for all $\alpha \geq \alpha_0$  and $q > 0$. 
(This is possible since $\cd(X)$ is finite.) Then
$\F \otimes \LL_\alpha$ is $0$-regular for $\alpha \geq \alpha_0$.
The conclusions then follow from the previous proposition.
\end{proof}

\section{Generalization of Serre Vanishing}\label{S:serre}

In this section we will prove Theorem~\ref{th:general-serre1}, 
a generalization of Serre's Vanishing Theorem.
It will allow us
to prove our desired implications in Theorem~\ref{th:maintheo}. 
Serre's Vanishing Theorem says that on a projective scheme with coherent sheaf $\F$ and ample $\LL$, 
the higher cohomology
of $\F \otimes \LL^{ m}$ vanishes for $m$ sufficiently large.
Our generalization says that for $m \geq m_0$, we also have the cohomological vanishing of
$\F \otimes \LL^{ m} \otimes \calN$, with $m_0$ \emph{independent} of $\calN$, where
$\calN$ is  numerically effective. This was proven for the case of
$X$ projective over an algebraically closed field in \cite[$\S$5]{Fuj2}
and we follow some of that proof. 

In essence, we will prove that a certain filter of invertible sheaves is an ample filter.
Thus we may use the results of $\S$\ref{S:ample1} to aid us. Let us
precisely define our filter of interest.

\begin{notation}
Let $X$ be a projective scheme over a noetherian ring $A$. Let
$\LL$ be an ample invertible sheaf  and let $\Lambda$ be a set 
of (isomorphism classes of) invertible sheaves on $X$. 
We will define a filter $(\LL, \Lambda)$ as follows.
As a set, $(\LL, \Lambda)$ is the collection of all
invertible sheaves $\LL^m \otimes \calN$ with $m \geq 0$ and $\calN \in \Lambda$.
For two elements $\calH_i$ of $(\LL, \Lambda)$, let $m_i$ 
be the maximum integer $m$ such that $\calH_i \cong \LL^{m} \otimes \calM$
for some $\calM \in \Lambda$. Then $\calH_1 < \calH_2$ if and only if $m_1 < m_2$.
This defines a partial ordering on $(\LL, \Lambda)$ which makes $(\LL, \Lambda)$
a filter of invertible sheaves.
\end{notation}

\begin{notation}\label{notation:Vq}
Let $A$ be a noetherian domain, let $\pi\colon X \to \Spec A$
be a projective morphism, let
$\LL$ be an ample invertible sheaf on $X$,  and let $\Lambda$ be a set
of (isomorphism classes of) invertible sheaves on $X$.
(For a locally closed subscheme $Y \subset X$, 
let $\Lambda\vert_Y = \{ \calN\vert_Y \colon \calN \in \Lambda \}$.)
Let $\F$ be a coherent sheaf on $X$.
For $q > 0$, 
we say that $V^q(\F, \LL, \Lambda)$ holds if there exists $m_0 = m(\F,q)$
and a non-empty open subscheme $U = U(\F, q) \subset \Spec A$ such that
$R^{q'}\pi_*( \F \otimes \LL^m \otimes \calN)\vert_V = 0$ for all
$m \geq m_0$, $q' \geq q$, $\calN \in \Lambda$, and all open subschemes $V \subset U$ such that 
$\calN\vert_{\pi^{-1}(V)}$ is $\pi\vert_{\pi^{-1}(V)}$-nef. 
If $\LL$ or $\Lambda$  are clear, we write $V^q(\F)$.
\end{notation}

In the vanishing conditions $V^q$ above, note that $m_0$ is independent of 
the particular open subscheme $V \subset U$ and that $\calN$
does not need to be nef over all of $\Spec A$. Thus the vanishing $V^q$
is in some sense stronger than the vanishing in Theorem~\ref{th:general-serre1}.
This will be necessary to reduce some of our work to the case of $A$ being
finitely generated over $\ZZ$.

We first prove reduction lemmas so we may work with schemes with the 
nice properties listed in Lemma~\ref{lem:geometrically-integral}.
It will be useful to replace $\Spec A$ with certain affine open subschemes.
This is possible because of \cite[Corollary~III.8.2]{Ha1}
\begin{equation}\label{eq:higher-direct-base-change}
R^q\pi_*(\F)\vert_U = R^q(\pi\vert_{\pi^{-1}(U)})_*(\F\vert_{\pi^{-1}(U)}).
\end{equation}
Many of our proofs implicitly use \eqref{eq:higher-direct-base-change}.
Our proofs also use descending induction on $q$, 
as we automatically have $V^q(\F)$ for any $\F$ and $q > \cd(X)$.

\begin{lemma}\label{lem:short-exact}
Let $X, A, \LL, \Lambda, q$ be as in Notation~\ref{notation:Vq}.
 Let
\[
0 \to \calK \to \F \to \G \to 0
\]
be a short exact sequence. Then for all $q > 0$,
\begin{enumerate}
\item $V^q(\calK, \LL, \Lambda)$ and $V^q(\G, \LL, \Lambda)$
imply $V^q(\F, \LL, \Lambda)$,
and we may assume $U(\F, q) = U(\calK, q) \cap U(\G, q)$,
\item $V^q(\F, \LL, \Lambda)$ and $V^{q+1}(\calK, \LL, \Lambda)$
imply $V^q(\G, \LL, \Lambda)$,
and we may assume $U(\G, q) = U(\F, q) \cap U(\calK, q+1)$, and
\item $V^q(\G, \LL, \Lambda)$ and $V^{q+1}(\F, \LL, \Lambda)$
imply $V^{q+1}(\calK, \LL, \Lambda)$,
and we may assume $U(\calK, q+1) = U(\G, q) \cap U(\F, q+1)$.
\end{enumerate}
\end{lemma}
\begin{proof}
The proof of each of the three statements is nearly the same.
For the first statement, set $U(\F, q) = U(\calK, q) \cap U(\G, q)$
and $m(\F, q) = \max\{ m(\calK, q), m(\G, q) \}$. Since $\Spec A$
is irreducible, the open set $U(\F, q)$ is non-empty.
\end{proof}

\begin{lemma}\label{lem:H1-enough}
Let $X, A, \LL, \Lambda, q$ be as in Notation~\ref{notation:Vq}.
Let $\F$ be a coherent sheaf on $X$. If $V^q(\OO_X, \LL, \Lambda)$ holds,
then $V^q(\F, \LL, \Lambda)$ holds
and one can take $U(\F, q) = U(\OO_X, q)$.
\end{lemma}
\begin{proof}
There is an exact sequence \eqref{eq:cokernel-of-ample}
\[
0 \to \calK \to \oplus_j \LL^{ -m} \to \F \to 0
\]
for some large $m$. By descending induction on $q$, we have that $V^{q+1}(\calK)$ holds
and $U(\calK, q+1) = U(\OO_X, q+1)$. Also, we may take $U(\OO_X, q+1) = U(\OO_X, q)$.
Now $V^q(\oplus_j \LL^{ -m})$ follows immediately from $V^q(\OO_X)$,
with $U(\oplus_j \LL^{ -m}, q) = U(\OO_X, q)$.
 We then have
$V^q(\F)$ by Lemma~\ref{lem:short-exact} and $U(\F, q) = U(\OO_X, q)$.
\end{proof}

We may now see a more direct connection between the vanishing $V^q$ and
ample filters as follows.

\begin{corollary}\label{cor:Vq-versus-ample-filter}
Let $X, A, \LL$ be as in Notation~\ref{notation:Vq}.
Let $\Lambda$ be a set of (not necessarily all) 
numerically effective invertible sheaves on $X$.
Then $V^1(\OO_X, \LL, \Lambda)$ holds if and only if
there exists an affine open subscheme $U\subset \Spec A$ such that
$(\LL\vert_{\pi^{-1}(U)}, \Lambda\vert_{\pi^{-1}(U)} )$
is an ample filter.
\end{corollary}
\begin{proof}
Suppose that $V^1(\OO_X, \LL, \Lambda)$ holds. We may take $U = U(\OO_X, 1)$ to
be affine. By Lemma~\ref{lem:H1-enough},
for any coherent sheaf $\F$ on $X$, we have $V^1(\F)$ and $U(\F, 1) = U$. 
For any $\calN \in \Lambda$, we have that $\calN\vert_{\pi^{-1}(U)}$
is nef, so the definition of $V^1$ gives the
vanishing necessary for $(\LL\vert_{\pi^{-1}(U)}, \Lambda\vert_{\pi^{-1}(U)} )$
to be an ample filter.

Now suppose that $(\LL\vert_{\pi^{-1}(U)}, \Lambda\vert_{\pi^{-1}(U)} )$
is an ample filter. By \eqref{eq:higher-direct-base-change}, 
we may replace $\Spec A$ with $U$
and $X$ with $\pi^{-1}(U)$. Then there exists $m_0$ such that
$R^q\pi_*(\LL^m \otimes \calN) = 0$ for all $q>0$, $m \geq m_0$ and $\calN \in \Lambda$.
If $V \subset U$ is an open subscheme, we trivially have the vanishing
necessary for $V^1(\OO_X, \LL, \Lambda)$.
\end{proof}

\begin{lemma}\label{lem:multiple-ample}
Let $X, A, \LL, \Lambda, q$ be as in Notation~\ref{notation:Vq}, 
and let $m > 0$. Then
$V^q(\OO_X, \LL, \Lambda)$ holds if and only if $V^q(\OO_X, \LL^m, \Lambda)$ holds.
\end{lemma}
\begin{proof}
By Lemma~\ref{lem:H1-enough}, this is actually a statement about the vanishing
$V^q(\F, \LL, \Lambda)$ and $V^q(\F, \LL^m, \Lambda)$
for all coherent sheaves $\F$. Given a coherent sheaf $\F$,
the statement $V^q(\F, \LL, \Lambda)$ obviously implies $V^q(\F, \LL^m, \Lambda)$.
Conversely, the statements $V^q(\F \otimes \LL^k, \LL^m, \Lambda)$, $k = 0, \dots, m-1$,
imply $V^q(\F, \LL, \Lambda)$, as in the proof of \cite[Proposition~II.7.5]{Ha1}.
\end{proof}

\begin{lemma}\label{lem:integral-reduction}
Let $X, A, \LL, \Lambda, q$ be as in Notation~\ref{notation:Vq}. If for all reduced, irreducible
components $X_i$ of $X$ the statement $V^q(\OO_{X_i}, \LL\vert_{X_i}, \Lambda\vert_{X_i})$ holds,
then $V^q(\OO_X, \LL, \Lambda)$ holds.
\end{lemma}
\begin{proof}
This is a standard reduction to the integral case, as in \cite[(5.10)]{Fuj2}
or \cite[Exercise~III.5.7]{Ha1}.
\end{proof}

\begin{lemma}\label{lem:cover-reduction}
Let $X, A, \LL, \Lambda, q$ be as in Notation~\ref{notation:Vq}, and let $X$ be integral. 
Let $\omega$ be a coherent sheaf on $X$ with $\Supp \omega = X$.
Suppose that $V^q(\omega,\LL, \Lambda)$ holds and 
for all coherent sheaves $\G$ with $\Supp \G \subsetneq X$,
the statement $V^q(\G,\LL, \Lambda)$ holds.
Then for all coherent sheaves $\F$, the statement $V^q(\F,\LL, \Lambda)$ holds.
\end{lemma}
\begin{proof}
By Lemma~\ref{lem:H1-enough}, it suffices to show $V^q(\OO_X)$ holds.
There exists $m$ sufficiently large so that $\calHom(\omega, \LL^m) \cong
\calHom(\omega, \OO_X) \otimes \LL^m$ is generated by global sections.
Now $\calHom(\omega, \OO_X) \neq 0$ since it is not zero at the generic point of $X$.
Thus there is a non-zero homomorphism $\phi\colon \omega \to \LL^m$.
Since $\LL^m$ is torsion-free, the sheaf $\G = \Coker(\phi)$ has $\Supp \G \subsetneq X$.

Consider the exact sequences
\begin{align}\label{eq:kernel-cokernel-image}
0 &\to \Ker(\phi) \to \omega \to \Image(\phi) \to 0 \\
\notag 0 &\to \Image(\phi) \to \LL^m \to \Coker(\phi) \to 0.
\end{align}
By descending induction on $q$, we may assume $V^{q+1}(\Ker(\phi))$ holds, thus
we have that $V^q(\Image(\phi))$ holds. Then
 $V^q(\LL^m)$ holds, since $V^q(\Coker(\phi))$ holds. Then obviously $V^q(\OO_X)$ holds.
\end{proof}

\begin{lemma}\label{lem:cover-reduction2}
Let $X, A, \LL, \Lambda, q$ be as in Notation~\ref{notation:Vq}. 
Let $f \colon X' \to X$ be a projective, surjective
morphism, let $\LL'$ be an ample invertible sheaf on $X'$, 
and let $\Lambda'$ be a set of invertible sheaves
on $X'$ such that 
\[
 \{ f^*\calN \colon \calN \in \Lambda \}
\cup \{ f^*\LL^m \colon m \geq 0 \} \subset \Lambda'.
\] 
If for all coherent sheaves $\F'$ on $X'$, the statement $V^q(\F', \LL', \Lambda')$
holds, then for all coherent sheaves $\F$ on $X$, the statement $V^q(\F, \LL, \Lambda)$
holds.  
\end{lemma}
\begin{proof}
If $Y \subset X$ is a closed  subscheme, then the hypotheses of the lemma are
satisfied by $X' \times_X Y \to Y$. So by Lemma~\ref{lem:integral-reduction}
and noetherian induction on $X$, we may assume that $X$ is integral and that
$V^q(\F)$ holds for all
coherent sheaves $\F$ on $X$ with $\Supp \F \subsetneq X$.

The invertible sheaf $\LL'$ is ample and hence $f$-ample, so take $m_0$ sufficiently
large so that $R^i f_*({\LL'}^m) = 0$ for $m \geq m_0$ and $i > 0$. 
Let $\pi' = \pi \circ f$.
Since $V^q(\OO_{X'}, \LL', \Lambda')$ holds, we may find an open subscheme $U \subset \Spec A$
and
 $m_1 \geq m_0$
so that $R^{q'}\pi'_*({\LL'}^m \otimes \calN')\vert_V = 0$ for all $m \geq m_1$, $q' \geq q$,
open $V \subset U$,
and $\calN' \in \Lambda'$ such that $\calN'\vert_{{\pi'}^{-1}(V)}$ 
is $\pi'\vert_{{\pi'}^{-1}(V)}$-nef.

If $\calN \in \Lambda$ and $\calN\vert_{\pi^{-1}(V)}$ is $\pi\vert_{\pi^{-1}(V)}$-nef, then
\[
f\vert_{{\pi'}^{-1}(V)}^*(\calN\vert_{\pi^{-1}(V)}) = (f^*\calN)\vert_{{\pi'}^{-1}(V)}
\]
is $\pi'\vert_{{\pi'}^{-1}(V)}$-nef by Lemma~\ref{lem:surjective-nef}.
So 
\[
R^{q'}\pi_*( f_*({\LL'}^{m_1}) \otimes \LL^m \otimes \calN)\vert_V
= R^{q'}\pi'_*( {\LL'}^{m_1} \otimes f^*\LL^m \otimes f^*\calN)\vert_V = 0
\]
for all $m \geq 0$, $q' \geq q$, open $V \subset U$ and $\calN \in \Lambda$
such that $\calN\vert_{\pi^{-1}(V)}$ is $\pi\vert_{\pi^{-1}(V)}$-nef.
Thus $V^q(f_*({\LL'}^{m_1}))$ holds. Now $\Supp f_*({\LL'}^{m_1}) = X$
since $f$ is surjective. So by Lemma~\ref{lem:cover-reduction}, we are done.
\end{proof}

We may now begin proving vanishing theorems akin to Theorem~\ref{th:general-serre1}.
We will first work with a ring $A$
which is of finite-type over $\ZZ$, then approximate a general noetherian ring $A$
with  finitely generated $\ZZ$-subalgebras.

\begin{proposition}\label{prop:general-serre-integer}
Let $A$ be a domain, finitely generated over $\ZZ$.
Let $X$ be projective over $A$, let $\LL$ be an ample invertible sheaf
on $X$, and let $\Lambda$ be the set of all invertible sheaves on $X$.
Then for all coherent sheaves $\F$ on $X$, the statement
$V^1(\F, \LL, \Lambda)$ holds.
\end{proposition}
\begin{proof}
By Lemma~\ref{lem:integral-reduction}, we may assume $X$ is integral.
If $\pi\colon X \to \Spec A$ is not surjective, then the claim
is trivial, as we can take $U(\F, 1)$ to be disjoint
from $\pi(X)$. So assume $\pi$ is surjective.
Let $d$ be the dimension of the generic fiber of $\pi$.
If $d = 0$, then $\pi$ is generically finite. In this case
we may replace $\Spec A$ with an affine open subscheme to assume $\pi$
is finite and hence $X$ is affine. Hence the proposition is trivial
\cite[Theorem~III.3.5]{Ha1}.
So assume that the proposition holds for all 
projective $Y \to \Spec A$ with generic fiber dimension $< d$.
Note that this implies that $V^q(\F)$ holds for all $\F$ with $\Supp \F \subsetneq X$.

By the definition of $V^q$, we may replace $S = \Spec A$ by an affine open
subscheme and then replace $X$ by a projective, surjective cover
by Lemma~\ref{lem:cover-reduction2}.
So if $A$ has characteristic $0$, we may assume $\pi$
is smooth and has geometrically integral fibers by 
Lemma~\ref{lem:geometrically-integral}.
Note that the generic fiber dimension does not change,
so we may still assume $V^1(\F)$ holds for all $\F$ with $\Supp \F \subsetneq X$.
Further, we may again replace $S$ with an affine open subscheme
to assume the morphism $S \to \Spec \ZZ$ is smooth 
\cite[Proposition~6.5]{Illusie-Frobenius} and 
assume that if $s \in S$ with $\character k(s) > 0$, then $\character k(s) > d$.

Let $\omega_{X/S} = \wedge^d \Omega^1_{X/S}$. The sheaf $\omega_{X/S}$ is invertible
since $\pi$ is smooth, and hence $\omega_{X/S}$ is flat over $S$. 
Let $V \subset S$ be an open subscheme, let $s \in V$
be closed in $V$, and let $X_s = X \times_S k(s)$ be the fiber. 
The residue field $k = k(s)$ is finite, hence perfect 
\cite[Proposition~6.4]{Illusie-Frobenius}. Let
$W_2(k)$ be the ring of second Witt vectors of $k$
\cite[3.9]{Illusie-Frobenius}.
 The closed immersion $\Spec k(s) \to V$
factors as $\Spec k \to \Spec W_2(k) \to V$,
since $V \to \Spec \ZZ$ is smooth \cite[2.2]{Illusie-Frobenius}.
Thus $X_s$ lifts to $W_2(k)$, i.e., there is scheme $X_1 = X \times_S W_2(k)$
with a smooth morphism
$X_1 \to \Spec W_2(k)$, such that $X_s = X_1 \times_{W_2(k)} k$.
So the Kodaira Vanishing Theorem holds for $X_s$ \cite[Theorem~5.8]{Illusie-Frobenius}
and thus
\[
H^q(X_s, \omega_{X_s/k(s)} \otimes \LL\vert_{X_s} \otimes \calN\vert_{X_s}) = 0
\]
for $q > 0$ and any $\calN$ such that $\calN\vert_{\pi^{-1}(V)}$ 
is $\pi\vert_{\pi^{-1}(V)}$-nef.
Now any open cover of the closed points of $V$ covers all of $V$ by
Lemma~\ref{lem:ZariskiSpace} and the invertible sheaf 
$\omega_{\pi^{-1}(V)/V}$ is flat over $V$, so 
\begin{multline*} 
R^q(\pi\vert_{\pi^{-1}(V)})_*(\omega_{\pi^{-1}(V)/V} \otimes \LL\vert_{\pi^{-1}(V)}
\otimes \calN\vert_{\pi^{-1}(V)}) \\ = R^q\pi_*(\omega_{X/S} \otimes \LL \otimes \calN)\vert_V = 0
\end{multline*}
for $q > 0$ by \cite[Theorem~III.12.11]{Ha1} and descending induction on $q$.
Thus we have that $V^1(\omega_{X/S}, \LL, \Lambda)$ holds.
So we are done with the characteristic $0$ case by Lemma~\ref{lem:cover-reduction}.

If $A$ has characteristic $p > 0$ (and hence  is finitely generated
over $K= \ZZ/p\ZZ$), then $X$ is quasi-projective over $K$,
so we may embed $X$ as an open subscheme of an integral scheme $\bar{X}$
which is projective over $K$.
Using
alteration of singularities \cite{deJong-Alterations}
and Lemma~\ref{lem:cover-reduction2},
we may assume $\bar{X}$ is a regular integral scheme, projective
and smooth
over $K$ (since $K$ is perfect).

Let $F\colon \bar{X} \to \bar{X}$ be the absolute Frobenius morphism.
Since $K$ is perfect, $F$ is a finite surjective morphism. Since $\bar{X}$ is regular,
$F$ is flat \cite[Corollary~V.3.6]{AltKle}, so there is a exact sequence of locally free sheaves
\[
0 \to \OO_{\bar{X}} \stackrel{\phi}{\to} F_*\OO_{\bar{X}}
\]
with the morphism $\phi$ locally given by $a \mapsto a^p$. The morphism $\phi$ 
remains injective
locally at $x \in \bar{X}$
upon tensoring with any residue field $k(x)$, so the cokernel of $\phi$ is also
locally free \cite[Exercise~II.5.8]{Ha1}.
Let $\omega_{\bar{X}} = \wedge^{\dim X} \Omega^1_{\bar{X}/K}$
and $\omega_X = \omega_{\bar{X}}\vert_X = \wedge^{\dim X} \Omega^1_{X/K}$
\cite[Theorem~VII.5.1]{AltKle}.
Dualizing via \cite[Exercises~III.6.10, 7.2]{Ha1}, there is a short exact sequence
of locally free sheaves
\[
0 \to \calK \to F_* \omega_{\bar{X}} \to \omega_{\bar{X}} \to 0.
\]
Restricting to $X$ gives
\[
0 \to \calK\vert_X \to F_* \omega_X \to \omega_X \to 0.
\]

By descending induction on $q$, assume
 $V^{q+1}(\calK\vert_X)$ and $V^{q+1}(\omega_X)$ hold.
 Let \[ U = U(\calK\vert_X, q+1) \cap U(\omega_X, q+1), \quad
 m_0 = \max \{ m(\omega_X, q+1), m(\calK\vert_X, q+1) \}. \]
 Then 
$R^{q'}\pi_*(\calK\vert_X \otimes \LL^m \otimes \calN)\vert_V = 0$ for all $m \geq m_0$,
$q' > q$, open subschemes $V \subset U$,
and $\calN \in \Lambda$ such that $\calN\vert_{\pi^{-1}(V)}$ is 
$\pi\vert_{\pi^{-1}(V)}$-nef. Then since $F$ is finite, 
\begin{multline*}
R^q\pi_*(\omega_X \otimes (\LL^m \otimes \calN)^{p^{e}})\vert_V
\cong R^q\pi_*( \omega_X \otimes F^*(\LL^m \otimes \calN)^{p^{e-1}})\vert_V \\
\cong R^q\pi_*( F_*\omega_X \otimes (\LL^m \otimes \calN)^{p^{e-1}})\vert_V
\to R^q\pi_*(\omega_X \otimes (\LL^m \otimes \calN)^{p^{e-1}})\vert_V \to 0
\end{multline*}
for all $e > 0$, $m \geq m_0$.
But the leftmost expression is $0$ for large $e$. Thus
$R^{q'}\pi_*( \omega_X \otimes \LL^m \otimes \calN)\vert_V = 0$ for $m \geq m_0$,
$q' \geq q$, open subschemes $V \subset U$, and $\calN \in \Lambda$ such that
$\calN\vert_{\pi^{-1}(V)}$ is 
$\pi\vert_{\pi^{-1}(V)}$-nef.
So $V^q(\omega_X)$ holds and we are again finished by Lemma~\ref{lem:cover-reduction}.
\end{proof}

\begin{remark}
As \cite[Lemma~5.8]{Fuj2} shows, it is possible to do the proof for $\character k = p > 0$
without assuming $\bar{X}$ is regular.
The cokernel of $F_* \omega_{\bar{X}} \to \omega_{\bar{X}}$ may not be zero,
but it is torsion, so the problem is solved by noetherian induction 
and exact sequences  similar to \eqref{eq:kernel-cokernel-image}.
In \cite{Fuj}, the regular case was presented to show
the main idea of the characteristic $p > 0$ proof.
However, it is interesting to note that now with alteration
of singularities, the proposition can
be reduced to the regular case in positive characteristic.
\end{remark}

We now move to the case of a general noetherian domain $A$. In the following proof,
we will see why it was necessary in Proposition~\ref{prop:general-serre-integer} to work
with $\Lambda$ equal to all invertible sheaves, instead of only numerically
effective invertible sheaves.

\begin{proposition}\label{prop:serre-finitely-generated}
Let $A$ be a noetherian domain, and let $X$ be projective over $A$.
Let $\LL, \calH_1, \dots, \calH_j$ be invertible sheaves on $X$
for some $j > 0$, and let $\LL$ be ample.
Let $\Lambda$ be the set of invertible sheaves on $X$ which are
numerically effective and which are in the subgroup of $\Pic X$
generated by the isomorphism classes of $\calH_i$, $i = 1, \dots, j$.
Then $V^1(\OO_X, \LL, \Lambda)$ holds.
\end{proposition}
\begin{proof}
The proposition is trivial if $\pi\colon X \to \Spec A$ is not surjective,
so assume this.

Let $S = \Spec A$.
There is a finitely generated $\ZZ$-subalgebra of $A$, call it $A_0$,
 a scheme $X_0$, a commutative diagram
\[
\xymatrix{ X_0 \times_{A_0} A \ar[r]^-{f} \ar[d]_\pi & X_0 \ar[d]^{\pi_0} \\
		S \ar[r]^g & 	S_0, }
\]
and invertible sheaves $\LL_0, \calH_{i,0}$ such that $S_0 = \Spec A_0$,
$X = X_0 \times_{A_0} A$, and $\LL \cong f^*\LL_0, \calH_i \cong f^*\calH_{i,0}$
 \cite[$\mathrm{IV}_3$, 8.9.1]{EGA}.
We may further assume that $\pi_0$ is projective and surjective
\cite[$\mathrm{IV}_3$, 8.10.5]{EGA} 
and that $\LL_0$ is ample \cite[$\mathrm{IV}_3$, 8.10.5.2]{EGA}.
By the definition of $V^1$, we may replace $S$ (and hence $S_0$)
with an affine open subscheme, so we may also assume $\pi_0$ (and hence $\pi$)
is flat \cite[$\mathrm{IV}_3$, 8.9.5]{EGA}.
Let $\Lambda_0$ be the set of all invertible sheaves on $X_0$.
By Proposition~\ref{prop:general-serre-integer},
we know $V^1(\OO_{X_0}, \LL_0, \Lambda_0)$ holds, so we may replace
$\Spec A_0$ with $U(\OO_{X_0}, 1)$.

Since all elements of $\Lambda$ are numerically effective, it suffices to show that
there exists $m_0$ which gives
the vanishing of $H^q(X, \LL^m \otimes \calN)$ for $q > 0$, $m \geq m_0$, 
and $\calN \in \Lambda$. Let $\calN \in \Lambda$, and let $\calN_0 \in \Lambda_0$
such that $\calN \cong f^*\calN_0$.
Let $s \in S$, and let $s_0 = g(s)$.
Then $(\LL \otimes \calN)\vert_{X_s}$ is ample,
and thus $(\LL_0 \otimes \calN_0)\vert_{X_0 \times_{A_0} k(s_0)}$ 
is ample by Lemma~\ref{lem:field-extension},
since $X_s = X \times_A k(s) = (X_0 \times_{A_0} k(s_0)) \times_{k(s_0)} k(s)$.
This is true for every $s_0 \in g(S)$, so there exists an open subscheme
$U \subset S_0$ such that $g(S) \subset U$ 
and $(\LL_0 \otimes \calN_0)\vert_{\pi_0^{-1}(U)}$
is $\pi_0\vert_{\pi_0^{-1}(U)}$-ample by Proposition~\ref{prop:fibers-ample}.

Let $m_0 = m(\OO_{X_0}, 1) + 1$, which does not depend on $\calN$.
Then we have $R^q(\pi_0)_*(\LL_0^m \otimes \calN_0)\vert_U = 0$ for $q > 0$ and $m \geq m_0$
by $V^1(\OO_{X_0})$. Since $\pi_0$ is flat, for each $s_0 \in g(S)$
we have $H^q(X_0 \times_{A_0} k(s_0), (\LL_0^m \otimes \calN_0)\vert_{X_0 \times_{A_0} k(s_0)}) = 0$
for $q > 0$, $m \geq m_0$, and $\calN \in \Lambda$ \cite[Theorem~III.12.11]{Ha1}.
The flat base change $\Spec k(s) \to \Spec k(s_0)$
gives $H^q(X \times_A k(s), (\LL^m \otimes \calN)\vert_{X \times_A k(s)})=0$
\cite[Proposition~III.9.3]{Ha1}. Another application of
\cite[Theorem~III.12.11]{Ha1}  gives the desired vanishing
of $H^q(X, \LL^m \otimes \calN)$.
\end{proof}

In view of Theorem~\ref{th:base1}, we have in some sense proven 
Theorem~\ref{th:general-serre1} ``up to numerical equivalence.''
We now prove a vanishing theorem for numerically trivial invertible sheaves.

\begin{proposition}\label{prop:serre-numerically-trivial}
Let $A$ be a noetherian domain, let $\pi\colon X \to \Spec A$ be a
flat, surjective, projective morphism with geometrically integral fibers, and 
let $\LL$ be a very ample invertible sheaf on $X$.
Let $\Lambda$ be the set of all numerically trivial invertible sheaves on $X$.
Then $( \LL, \Lambda)$ is an ample filter.
\end{proposition}
\begin{proof}
Let $\calN$ be a numerically trivial sheaf on $X$, and let $s \in \Spec A$.
Using the notation of Remark~\ref{remark:algebraic-closure-fiber},
we see that $f(m) = \chi( \LL_{\bar{s}}^m \otimes \calN_{\bar{s}})$ is a polynomial in $m$. 
Note that $f(m)$ does not depend on $\calN$ since $\calN_{\bar{s}}$ is
numerically trivial.
By \cite[p.~312, Theorem~3]{K}, there exists $m_0$, which depends only
on the coefficients of $f(m)$, such that
\[
H^q(X_{\bar{s}},  \LL_{\bar{s}}^m \otimes \calN_{\bar{s}}) = 0
\]
for $m \geq m_0$ and $q > 0$. 

But $f(m)$ does not depend
on $s$, as noted in Remark~\ref{remark:algebraic-closure-fiber}.
Then using the flat base change $\Spec \ov{k(s)} \to \Spec k(s)$
and \cite[Theorem~III.12.11]{Ha1}, we have that
\[
H^q(X,  \LL^m \otimes \calN ) = 0
\]
for $m \geq m_0$, $q > 0$, and all numerically trivial $\calN$.
Thus by Proposition~\ref{prop:vanishing-ample-powers-enough}, $(\LL, \Lambda)$
is an ample filter.
\end{proof}

Now via the following lemma, we tie together the vanishing 
in Proposition~\ref{prop:serre-finitely-generated}
with that in Proposition~\ref{prop:serre-numerically-trivial}.

\begin{lemma}\label{lem:ugly-lemma}
Let $k$ be a field, and let $X$ be an equidimensional scheme, projective
over $k$. Let $\Lambda$ be a set of (not necessarily all)
numerically trivial invertible sheaves on $X$, containing $\OO_X$
and closed under inverses. 
Let $\LL$ be an invertible sheaf on $X$
such that $\LL \otimes \calN$ is very ample for all 
$\calN \in \Lambda$.
Let $V$ be a closed subscheme of $X$, let $t > 0$, and let $r(V,t) = 2^{\dim V - t+1}-1$
if $\dim V \geq t-1$, and let $r(V, t) = 0$ otherwise.
Let $\calH$ be an invertible sheaf on $X$
such that 
\[
H^{q'}(V, \OO_V \otimes \calH \otimes \LL^{m-1})= 0
\]
 for $q' \geq q$
and $m \geq m_0$.
Then 
\[
H^{q'}(V, \OO_V \otimes \calH \otimes \calN \otimes \LL^{m-1+r(V,q')})= 0
\]
for
$q' \geq q$, $m \geq m_0$, and
all $\calN \in \Lambda$.
\end{lemma}
\begin{proof}
We proceed by induction on $\dim V$, the claim
being trivial when $\dim V = 0$.
Let $q > 0$. We may assume $q \leq \dim V$.
By descending induction on $q$, we may also assume
$H^{q'}(V, \OO_V \otimes \calH \otimes \calN \otimes \LL^{m-1+r(V,q')}) = 0$
for $q' > q$, $m \geq m_0$, and  $\calN \in \Lambda$.

Since $\LL \otimes \calN$ is very ample, there is an effective Cartier divisor $D \subset V$
with an exact sequence
\begin{equation}\label{eq:ugly-sequence}
0 \to \OO_V \otimes \calH \otimes \calN^{-1} \otimes \LL^{m-1} 
\to \OO_V \otimes \calH  \otimes \LL^{m}
\to \OO_D \otimes \calH  \otimes \LL^{m}
\to 0.
\end{equation}
Tensoring \eqref{eq:ugly-sequence} with $\LL^{r(V, q+1)}$
and examining the related long exact sequence,
we have $H^{q'}(D, \OO_D \otimes \calH \otimes \LL^{m+r(V,q+1)}) = 0$
for $q' \geq q$ and $m \geq m_0$, since $\calN^{-1} \in \Lambda$.

Then by induction on $\dim V$, we have
\[
H^q(D, \OO_D \otimes \calH \otimes \LL^{m+r(V, q+1) + r(D, q)} \otimes \calM) = 0
\]
for $m \geq m_0$ and any $\calM \in \Lambda$.
Now $\dim D = \dim V - 1$ and
\begin{multline*} 
r(V, q+ 1) + r(D, q) = 2^{\dim V - q} - 1 + 2^{\dim D - q + 1} - 1 \\
= 2^{\dim V - q + 1} - 2 = r(V, q) - 1.
\end{multline*}
Note that since $\dim V \geq q$, we have $r(V, q) -1 \geq 0$.
Tensoring \eqref{eq:ugly-sequence} with $\LL^{r(V, q) - 1} \otimes \calN$,
we have $H^q(V, \OO_V \otimes \calH \otimes \LL^{m-1+ r(V, q)} \otimes \calN) = 0$
for $m \geq m_0$. Since this holds for any $\calN \in \Lambda$,
we have proven the lemma.
\end{proof}

We may now finally prove Theorem~\ref{th:general-serre1}.

\begin{proof}
Let $S = \Spec A$, and
let $\Lambda$ be the set of nef invertible sheaves on $X$. We wish
to show that $(\LL, \Lambda)$ is an ample filter. By 
Proposition~\ref{prop:filter-reductions}, we may assume $X$ is an integral
scheme. Let $\pi\colon X \to S$ be the structure morphism. We may
replace $S$ with $\pi(X)$ and hence assume that $A$ is a domain
and that $\pi$ is surjective.

The theorem is trivial if $X = \emptyset$. 
Note that for a closed subscheme $Y$ of $X$, that
$\Lambda\vert_Y$ is a subset of the set of all numerically effective sheaves on $Y$.
So by noetherian induction on $X$,
we may assume that for any proper closed subscheme $Y$ of $X$ (in particular
$Y = \pi^{-1}(T)$ for some proper closed subscheme $T \subset S$),
that $(\LL\vert_Y, \Lambda\vert_Y)$ is an ample filter. So by
Proposition~\ref{prop:ample-cover-implies-ample-filter}, we need only show
that there exists an affine open subscheme $U \subset S$ such that
$(\LL\vert_{\pi^{-1}(U)}, \Lambda\vert_{\pi^{-1}(U)})$ is an ample filter.
Or equivalently by Corollary~\ref{cor:Vq-versus-ample-filter},
we need to show that $V^1(\OO_X, \LL, \Lambda)$ holds.

By the definition of $V^1$, we may replace $S$ with any 
affine open subscheme, and by Lemmas~\ref{lem:cover-reduction2}
and \ref{lem:geometrically-integral}, we may assume $\pi$
is flat and has geometrically integral fibers. Let $d$
be the dimension of the generic fiber of $\pi$.

Let $\Theta$ be the set of all numerically trivial invertible sheaves
on $X$. By Proposition~\ref{prop:serre-numerically-trivial}, $(\LL, \Theta)$ 
is an ample filter.
So by Corollary~\ref{cor:general-serre-global-sections}, there exists $m$
such that $\LL^m \otimes \calM$ is very ample
for all $\calM \in \Theta$. By Lemma~\ref{lem:multiple-ample}, we replace $\LL$
by $\LL^m$.

Choose (finitely many) $\ZZ$-generators $\calH_i$, $i = 1,\dots, \rho(X)$, 
of $A^1(X)$ and let
$\Sigma$ be the set of
all nef invertible sheaves 
in the subgroup of $\Pic X$ generated by $\calH_i$, $i = 1,\dots, \rho(X)$.
Then by Proposition~\ref{prop:serre-finitely-generated}, we may find
an affine open subscheme $U \subset S$ and $m_0$ such that
\[
R^q \pi_*(\LL^m \otimes \calH)\vert_U = 0
\]
for $q> 0$, $m \geq m_0$, and $\calH \in \Sigma$.
Since $\pi$ is flat, for $s \in U$ we then have \cite[Theorem~III.12.11]{Ha1}
\[
H^q(X_s, \LL_s^m \otimes \calH_s) = 0.
\]
Then by Lemma~\ref{lem:ugly-lemma}, we have
\[
H^q(X_s, \LL_s^{m+2^d} \otimes \calH_s \otimes \calM_s) = 0
\]
for $q> 0$, $m \geq m_0$, $\calH \in \Sigma$, and $\calM \in \Theta$.
Now any $\calN \in \Lambda$ can be written as $\calN \cong \calH \otimes \calM$
for some $\calH \in \Sigma$, $\calM \in \Theta$.
So another application of \cite[Theorem~III.12.11]{Ha1}
gives
\[
R^q\pi_*(\LL^m \otimes \calN)\vert_U = 0
\]
for $q> 0$, $m \geq m_0 + 2^d$, and $\calN \in \Lambda$.
Thus $V^1(\OO_X, \LL, \Lambda)$ holds, as desired.
\end{proof}

Theorem~\ref{th:general-serre1} is the best possible in the following sense.

\begin{proposition}\label{prop:best-set}
Let $A$ be a noetherian ring,
let $X$ be projective  over $A$, and let
$\LL$ be an ample
invertible sheaf on $X$. Let $\Lambda$ be a set of invertible sheaves
on $X$ such that if $\calN \in \Lambda$, then $\calN^j \in \Lambda$ for
all $j > 0$. Suppose that for all coherent sheaves $\F$, there exists $m_0$
such that
\[
H^q(X, \F \otimes \LL^m \otimes \calN) = 0
\]
for all $m \geq m_0$, $q > 0$, and all $\calN \in \Lambda$.
Then $\calN$ is numerically effective for all $\calN \in \Lambda$.
\end{proposition}
\begin{proof}
Suppose that $\calN \in \Lambda$
is not numerically effective. Then there exists a contracted 
integral curve $f\colon C \hookrightarrow X$
with  $(\calN.C) < 0$. Then by the Riemann-Roch formula, for any fixed $m$,
we can choose $j$ sufficiently large so that
\[
H^1(C, \OO_C \otimes f^*(\LL^{ m} \otimes \calN^{ j})) \neq 0.
\]
But $\calN^j \in \Lambda$,
so this is a contradiction. Therefore, $\calN$ is numerically effective.
\end{proof}

We also have an immediate useful corollary
to Theorem~\ref{th:general-serre1}, 
via Corollary~\ref{cor:general-serre-global-sections}.

\begin{corollary}\label{cor:general-serre-global-sections2}
Let $A$ be a  noetherian ring,
let $X$ be a projective scheme over $A$,  let $\LL$
be an ample
invertible sheaf on $X$ with $\LL^n$ very ample,
 and let $\F$ be a coherent sheaf on $X$.
There exists $m_0$ such that
\begin{enumerate}
\item The natural map 
$H^0(X, \F \otimes \LL^m \otimes \calN) \otimes H^0(X, \LL^n)
\to H^0(X, \F \otimes \LL^{m+n} \otimes \calN)$ is surjective and
\item $\F \otimes \LL^m \otimes \calN$ is generated by global sections
\end{enumerate}
for $m \geq m_0$ and all numerically effective invertible sheaves $\calN$.\qed
\end{corollary}

\section{Ample filters II}\label{S:ample-posets}

In this section we will prove Theorem~\ref{th:maintheo}.
In \cite[Proposition~3.2]{AV}, it was first 
shown that  \eqref{th:main1} implied \eqref{th:main3}, and
 \eqref{th:main3} implied  \eqref{th:main4}, in the case of an ample sequence
and
$X$ projective over a  field $k$. The proof makes strong use of the projectivity
of $X$ and also requires the vanishing of cohomology in  \eqref{th:main1}; the surjective
map of global sections in  \eqref{th:main2} would not have been strong enough
for this method. That \eqref{th:main4} implies \eqref{th:main1} was
first noted in \cite[Proposition~2.3]{Ke1}, in the case of an algebraically closed field
and a certain ample sequence.

We begin our proof of Theorem~\ref{th:maintheo} by
using the results of $\S$\ref{S:serre}.

\begin{proposition}\label{prop:A4-implies-A1-A3}
Let $X$ be a scheme, proper over a noetherian ring $A$. 
 Let $\{ \LL_\alpha \}$
be a filter of invertible sheaves. If  $\{ \LL_\alpha \}$ satisfies {\eqref{th:main4}}, 
then $\{ \LL_\alpha \}$ satisfies {\eqref{th:main1}} and {\eqref{th:main3}}.
\end{proposition}
\begin{proof}
Let $\calH$ be an ample invertible sheaf on $X$, and
let $\F$ be a coherent sheaf on $X$. By Theorem~\ref{th:general-serre1} and
Corollary~\ref{cor:general-serre-global-sections2}, there exists $m$ such that
$\F \otimes \calH^m \otimes \calN$ has vanishing higher cohomology and is
generated by global sections for any nef invertible sheaf $\calN$.
By \eqref{th:main4}, there exists $\alpha_0$ such that $\calH^{-m} \otimes \LL_\alpha$
is ample (hence nef) for $\alpha \geq \alpha_0$. Thus $\F \otimes \LL_\alpha$
has the desired properties for $\alpha \geq \alpha_0$.
\end{proof}

The statements \eqref{th:main1}--\eqref{th:main3} do not immediately imply
that $X$ is projective over $A$, so we may not assume that $X$ has ample
invertible sheaves. It will be much easier to argue certain sheaves are nef.
Thus we need to prove that a fifth statement is equivalent to
statement~\eqref{th:main4} of  Theorem~\ref{th:maintheo}. This is a natural generalization
of Kleiman's criterion for ampleness.

\begin{proposition}\label{prop:amplenef}
Let $X$ be a scheme, proper over a  noetherian ring $A$, and let $\{\LL_\alpha \}$ be
a filter of invertible sheaves. Then the following are equivalent: 
\begin{enumerate}
\setcounter{enumi}{3}
\renewcommand{\theenumi}{A\arabic{enumi}}
\item\label{prop:amplenef1} For all invertible sheaves $\calH$, there exists
$\alpha_0$ such that $\calH \otimes \LL_\alpha$ is an ample invertible sheaf
for $\alpha \geq \alpha_0$.
\item\label{prop:amplenef2} The scheme $X$ is  quasi-divisorial, and for all 
invertible sheaves $\calH$, there exists
$\alpha_0$ such that $\calH \otimes \LL_\alpha$ is a numerically effective invertible sheaf
for $\alpha \geq \alpha_0$.
\end{enumerate}
\end{proposition}
\begin{proof} \eqref{prop:amplenef1} $\implies$ \eqref{prop:amplenef2} is clear
since $X$ is projective as it has an ample
invertible sheaf \cite[Remark~II.5.16.1]{Ha1}. 
Hence $X$ is  quasi-divisorial.
Further, any ample invertible sheaf
is necessarily numerically effective.
Thus \eqref{prop:amplenef2} holds.

Now assume \eqref{prop:amplenef2}. 
Let $A^1 = A^1(X)$.
Then $A^1$ is a finitely generated free abelian group by Theorem~\ref{th:base},
with rank $\rho(X)$.
Let $\calH_1, \dots, \calH_{\rho(X)}$ be a $\ZZ$-basis for $A^1$. 
For $v \in A^1 \otimes \RR$, we may write
\[
v = a_1 \calH_1 + a_2 \calH_1^{-1} + \dots + a_{2\rho(X)} \calH_{\rho(X)}^{-1}
\]
with $a_i \geq 0$ in the additive notation of $A^1 \otimes \RR$.

We choose  an arbitrary invertible sheaf $\calH$.
Since $\{ \LL_\alpha \}$ is a filter, we may choose $\alpha_0$ large enough
so that for $i = 1, \dots, \rho(X)$ and $\alpha \geq \alpha_0$, the invertible sheaves 
$\calH \otimes \LL_\alpha$ and
$\calH_i^{\pm 1} \otimes \calH \otimes \LL_\alpha$
are numerically effective. Set $a = \sum a_i$. Then
\[
v + a (\calH \otimes \LL_\alpha) + b (\calH \otimes \LL_\alpha)
\]
is in $K$, the cone generated by numerically effective invertible sheaves, for 
$\alpha \geq \alpha_0$ and $b\geq 0$. Thus by Lemma~\ref{lem:interior},
we have $\calH \otimes \LL_\alpha \in \interior K$ and hence $\calH \otimes \LL_\alpha$
is ample by Theorem~\ref{th:Kleiman}, as desired.
\end{proof}

\begin{lemma}\label{lem:finitebasenef}
Let $X$ be proper over a field $k$, and let $\LL$ be an invertible sheaf on $X$. 
Suppose that the base locus $B$ of $\LL$, the points at which
$\LL$ is not generated by global sections, is zero-dimensional or empty. 
Then $\LL$ is numerically effective.
\end{lemma}
\begin{proof}
Let $C \subset X$ be an integral curve. Since $B$ is zero-dimensional or empty, 
$C$ is not a subset of $B$.
Thus there exists $t \in H^0(X, \LL)$ such that $U_t \cap C \neq \emptyset$, where
$U_t$ is the (open) set of points $x \in X$ such that the stalk $t_x$ of $t$ at $x$ is 
not contained in $\scrm_x \LL_x$, where $\scrm_x$ is
the maximal ideal of $\OO_{X,x}$ \cite[Lemma~II.5.14]{Ha1}.
Thus $(\LL.C) \geq 0$  \cite[Example~12.1.2]{Fu}.
\end{proof}

\begin{lemma}\label{lem:affinebasenef}
Let $X$ be proper over a noetherian ring $A$, and let $\LL$ be an invertible sheaf on $X$.
Suppose that there exists an affine open subscheme $U$ of $X$ such that
$\LL$ is generated by global sections at all $x \in X \setminus U$.
Then $\LL$ is numerically effective.
\end{lemma}
\begin{proof} Let $V = X \setminus U$ with the reduced induced closed
subscheme structure. We choose a closed point $s \in \Spec A$. 
Then $X_s = X \times_A k(s)$ is a closed subscheme of $X$. For all
$x \in X_s$ there is a commutative diagram
\[
\xymatrix{ \LL \ar[r] \ar[d] & \OO_{X_s} \otimes \LL \ar[d] \\
			\LL_x \ar[r] & (\OO_{X_s} \otimes \LL)_x \ar[r] & 0. }
\]
So $\OO_{X_s} \otimes \LL$ is generated by global sections for
all $x \in V \times_A k(s)$. Thus the base locus $B_s$ of points
at which $\OO_{X_s} \otimes \LL$ is not generated by global sections
is contained in the open affine subscheme $U \times_A k(s)$.
So $B_s$ is an affine scheme since $B_s$ is a closed subscheme of $U \times_A k(s)$.
But since $B_s$ is a closed subscheme of $X_s$, the natural morphism
$B_s \to \Spec k(s)$ is affine and proper, hence finite \cite[Exercise~II.4.6]{Ha1}.
So $B_s$ is zero-dimensional or empty, and so $\LL\vert_{X_s}$ is numerically effective
by Lemma~\ref{lem:finitebasenef}. Since this is true for every closed point
$s \in S$, we must have that $\LL$ is numerically effective.
\end{proof}

\begin{proposition}\label{prop:A3-implies-A4}
Let $X$ be a scheme, proper over a noetherian ring $A$. 
 Let $\{ \LL_\alpha \}$
be a filter of invertible sheaves. If  $\{ \LL_\alpha \}$ satisfies {\eqref{th:main3}}, 
then $\{ \LL_\alpha \}$ satisfies {\eqref{th:main4}}. 
\end{proposition}
\begin{proof} 
We show \eqref{prop:amplenef2} holds and 
thus \eqref{th:main4} holds by Proposition~\ref{prop:amplenef}.
First, we must show that $X$ is quasi-divisorial.
Let $V$ be a closed, contracted, integral
subscheme of $X$. There exists $\alpha_0$ such that
$\LL_\alpha\vert_V$ is generated by global sections for
$\alpha \geq \alpha_0$. 
Since $V$ is an integral scheme, each $\LL_\alpha\vert_V$ must
define an effective Cartier divisor \cite[Proposition~II.6.15]{Ha1}.
Suppose that all these Cartier divisors are zero.
Then $\LL_\alpha\vert_V \cong \OO_V$ for all $\alpha \geq \alpha_0$
and thus
all coherent sheaves $\F$ on $V$ are generated by global sections. But then
$V$ is affine and proper over $\Spec k(s)$. Hence $V$ is zero-dimensional
 \cite[Exercise~II.4.6]{Ha1} and hence $V$ is an integral point. 
So $X$ is 
quasi-divisorial.

Given any invertible sheaf $\calH$ on $X$, there is $\alpha_1$
such that $\calH \otimes \LL_\alpha$ is generated by global sections
for $\alpha \geq \alpha_1$. A trivial application of Lemma~\ref{lem:affinebasenef}
shows that $\calH \otimes \LL_\alpha$ is numerically effective, as desired.
\end{proof}

\begin{proposition}\label{prop:A2-implies-A3}
Let $X$ be a scheme, proper over a noetherian ring $A$. 
 Let $\{ \LL_\alpha \}$
be a filter of invertible sheaves. If  $\{ \LL_\alpha \}$ satisfies {\eqref{th:main2}}, 
then $\{ \LL_\alpha \}$ satisfies  {\eqref{th:main3}}.
\end{proposition}
\begin{proof}
Let $\F$ be a coherent sheaf. If $\Supp \F = \emptyset$, then the claim is obvious.
So by noetherian induction on $X$, suppose that the proposition holds on all proper closed
subschemes $V$ of $X$.

Since \eqref{th:main3} is equivalent to \eqref{prop:amplenef2}
(by Propositions~\ref{prop:A4-implies-A1-A3}, \ref{prop:amplenef}, and \ref{prop:A3-implies-A4}),
we may show that \eqref{prop:amplenef2} holds.
Let $W$ be a closed, contracted, integral
subscheme of $X$ with $\pi(W) = s$, where $\pi$ is the structure morphism.
We choose a closed point
$x \in W$.
Then there exists $\alpha_0$ such that
\[
H^0(W,\LL_\alpha\vert_W) \to H^0(\{ x \},  \LL_\alpha\vert_{\{ x \} })
\]
is surjective for $\alpha \geq \alpha_0$, where $\{ x \}$ is the 
 reduced closed subscheme defined by the point $x$.  
So $\dim_{k(s)} H^0(W, \LL_\alpha\vert_W) \geq 1$.
Since $W$ is an integral scheme, each $\LL_\alpha\vert_W$ must
define an effective Cartier divisor \cite[Proposition~II.6.15]{Ha1}.
Suppose that all these Cartier divisors are zero.
Then $\LL_\alpha\vert_W \cong \OO_W$ for all $\alpha \geq \alpha_0$ and
one can argue that $\OO_W$ is an ample invertible sheaf on $W$,
using the proof of \cite[Proposition~III.5.3]{Ha1}. But then
$W$ is affine and proper over $\Spec k(s)$. Hence $W$ is zero-dimensional 
\cite[Exercise~II.4.6]{Ha1} and hence $W$ is an integral point.
So $X$ is 
quasi-divisorial.

Now let $U$ be an affine open subscheme of $X$. Set $V = X \setminus U$
with the reduced induced subscheme structure.
Let $\calH$ be an invertible sheaf on $X$.
Since \eqref{th:main3} holds for
$\{ \LL_\alpha\vert_V \}$ there is $\alpha_1$ such that
$(\calH \otimes \LL_\alpha)\vert_V$ 
is generated by global sections for $\alpha \geq \alpha_1$.

There exists $\alpha_2 \geq \alpha_1$ such that 
$H^0(X, \calH \otimes \LL_\alpha) \to H^0(V, (\calH \otimes \LL_\alpha)\vert_V)$ is
an epimorphism for $\alpha \geq \alpha_2$. Thus if $x \in V$, then the stalk
$\calH \otimes \LL_\alpha \otimes \OO_{V,x}$ is generated by $H^0(X, \calH \otimes \LL_\alpha)$.
Nakayama's Lemma implies $H^0(X, \calH \otimes \LL_\alpha)$ also generates the stalk
$\calH \otimes \LL_\alpha \otimes \OO_{X,x}$.
So by Lemma~\ref{lem:affinebasenef}, $\calH \otimes \LL_\alpha$ is
numerically effective for $\alpha \geq \alpha_2$, as we wished to show.
Thus \eqref{prop:amplenef2} holds for $\{ \LL_\alpha \}$.
\end{proof}

We may now summarize the proof of Theorem~\ref{th:maintheo}.

\begin{proof}
\eqref{th:main1} $\implies$ \eqref{th:main2}:
Apply the first statement to get the vanishing of
\[
H^1(X, \Ker(\F \twoheadrightarrow \G) \otimes \LL_\alpha).
\]
The desired surjectivity then follows from the natural long exact sequence.

\eqref{th:main2} $\implies$ \eqref{th:main3}:
This is Proposition~\ref{prop:A2-implies-A3}.

\eqref{th:main3} $\implies$ \eqref{th:main4}:
This is Proposition~\ref{prop:A3-implies-A4}.

\eqref{th:main4} $\implies$ \eqref{th:main1}:
This is Proposition~\ref{prop:A4-implies-A1-A3}.
\end{proof}

\begin{corollary}\label{cor:projective}
Let $X$ be proper over a noetherian ring $A$. Then $X$ is projective
over $A$ if and only if $X$ has an ample filter of invertible sheaves.
\end{corollary}
\begin{proof}
If $X$ is projective, it has an ample invertible sheaf $\LL$ and
the filter $\{ \LL, \LL^{ 2}, \dots \}$ is an ample filter. Conversely, if $X$ has
an ample filter, it has an ample invertible sheaf by Theorem~\ref{th:maintheo}.
Thus $X$ is projective \cite[Remark~II.5.16.1]{Ha1}.
\end{proof}

\begin{remark}\label{remark:divisorial}
Let $X$ be a separated noetherian scheme.
If $X$ is covered by affine open complements of Cartier divisors, then
$X$ is called \emph{divisorial}. All projective schemes and all regular proper
integral schemes
(over an affine base) are
divisorial \cite[Pf. of Theorem~VI.2.19]{KolRat}, and any divisorial scheme
is quasi-divisorial.
The above corollary is not true in the divisorial case 
if one replaces ``ample filter of invertible sheaves''
with ``ample filter of non-zero torsion-free coherent subsheaves of invertible sheaves.''
To be more exact, let $X$ be a normal, divisorial, proper scheme over an algebraically
closed field. Then there exists an invertible sheaf $\LL$ and a sequence of
non-zero coherent sheaves of ideals $\{ \calI_m \}$ such that for all coherent sheaves $\F$
on $X$, there exists $m_0$ such that
\[
H^q(X, \F \otimes \calI_m \otimes \LL^{ m}) = 0
\]
for all $q > 0$ and $m \geq m_0$ \cite[Corollary~6]{B-CohDiv}.
Recently, it was shown that if $X$ is divisorial and proper over a noetherian ring $A$,
then there exists a non-zero subsheaf $\calK$ of an invertible sheaf $\calH$ such that
for all coherent $\F$, there exists $m_0$ such that $H^q(X, \F \otimes \calK^m) = 0$
for all $q> 0$ and $m \geq m_0$ \cite[Theorem~5.3]{BS-AmpleFam}.
\end{remark}

To conclude this section, we examine other related conditions on
a filter of invertible sheaves.

\begin{proposition}\label{prop:ample-poset-v-ample}
Let $X$ be a scheme, proper over a noetherian ring $A$. Let $\{ \LL_\alpha \}$
be a filter of invertible sheaves. The filter $\{ \LL_\alpha \}$
is an ample filter if and only if
for all invertible sheaves $\calH$, there exists $\alpha_0$ such that
	$\calH \otimes \LL_\alpha$ is a very ample invertible sheaf 
	for $\Spec A$ when $\alpha \geq \alpha_0$.
\end{proposition}
\begin{proof}
Suppose that $\{ \LL_\alpha \}$ is an ample filter. By Corollary~\ref{cor:projective},
$X$ is projective over $A$ and hence has a very ample invertible sheaf $\OO_X(1)$.
Let $\calH$ be an invertible sheaf.
For $\alpha$ sufficiently large, $\OO_X(-1) \otimes \calH \otimes \LL_\alpha$ is generated
by global sections. Thus 
\[
\OO_X(1) \otimes \OO_X(-1) \otimes \calH \otimes \LL_\alpha \cong \calH \otimes \LL_\alpha
\]
is very ample \cite[Exercise~II.7.5]{Ha1}. 

The converse is clear, since any very ample invertible sheaf is ample, so
\eqref{th:main4} holds for $\{ \LL_\alpha \}$.
\end{proof}

These propositions regarding ample filters also have a relative form. The
proofs are as in \cite[Theorem~III.8.8]{Ha1}. We state
a relative form of Theorem~\ref{th:maintheo}.

\begin{theorem}
Let $S$ be a noetherian scheme,
let $\pi\colon X \to S$ be a proper morphism, and
let $\{ \LL_\alpha \}$
be a filter of invertible sheaves on $X$. 
Then the following are equivalent:
\begin{enumerate}
	\item For all coherent sheaves $\F$ on $X$, there exists $\alpha_0$ such that
	$R^q\pi_*(\F \otimes \LL_\alpha)=0$ for all $q > 0$ and $\alpha \geq \alpha_0$,
	i.e., $\{ \LL_\alpha \}$ is a $\pi$-ample filter.
	\item For all coherent sheaves $\F, \G$ on $X$ with epimorphism 
		$\F \twoheadrightarrow \G$, there exists $\alpha_0$ such that the natural map
		\[
			\pi_*(\F \otimes \LL_\alpha) \to \pi_*(\G \otimes \LL_\alpha)
		\]
		is an epimorphism for $\alpha \geq \alpha_0$.
	\item For all coherent sheaves $\F$, there exists $\alpha_0$ such that
	the natural morphism
		$\pi^*\pi_*(\F \otimes \LL_\alpha) \to \F \otimes \LL_\alpha$ 
		is an epimorphism for $\alpha \geq \alpha_0$.
	\item For all invertible sheaves $\calH$ on $X$, there exists $\alpha_0$ such that
	$\calH \otimes \LL_\alpha$ is a $\pi$-ample invertible sheaf for $\alpha \geq \alpha_0$.\qed
\end{enumerate}
\end{theorem}

\section{Twisted homogeneous coordinate rings}\label{S:rings}

Throughout this section, the scheme $X$ will be proper over a commutative
noetherian ring $A$.
We now generalize the results of \cite{Ke1} to this case. Specifically,
we show that left and right $\gs$-ampleness are equivalent in this
case and thus the associated twisted homogeneous coordinate ring
is left and right noetherian.

 First we must briefly review
the concept of a noncommutative projective scheme from \cite{AZ}.
Let $R$ be a noncommutative, right noetherian, $\NN$-graded ring, let $\gr R$ be the category
of finitely generated, graded right $R$-modules, let $\tors R$ be the full subcategory
of torsion submodules, and let $\qgr R$ be the quotient category
$\gr R/\tors R$. Let $\pi\colon \gr R \to \qgr R$ be the quotient functor.
Then the pair $\proj R = (\qgr R, \pi R)$ is said to be a \emph{noncommutative
projective scheme}. We will work with rings $R$ such that
$\proj R$ is  equivalent to $(\coh(X), \OO_X)$, where $\coh(X)$ is
the category of coherent sheaves on $X$. By saying that $\proj R \cong (\coh(X), \OO_X)$
we mean that there is a category equivalence $f\colon \qgr R \to \coh(X)$ and
$f(\pi R) \cong \OO_X$.

Given an $A$-linear abelian category $\calC$, arbitrary object $\OO$, and autoequivalence
$s$, one can define a homogeneous coordinate ring
\[
R = \Gamma(\calC, \OO, s)_{\geq 0} = \bigoplus_{i=0}^\infty \Hom(\OO, s^i\OO)
\]
with multiplication given by composition of homomorphisms. That is,
given $a \in R_m, b \in R_n$,  we have $a\cdot b = s^n(a) \circ b \in \Hom(\OO, s^{n+m}\OO)$.

Without loss of generality, one may assume $s$ is a category automorphism, i.e., there
is an inverse autoequivalence $s^{-1}$ \cite[Proposition~4.2]{AZ}. We then have a concept of
ample autoequivalence.

\begin{definition}\label{def:ample-auto}
Let $\calC$ be an $A$-linear abelian category. A pair $(\OO, s)$, with $\OO \in \calC$
and  an autoequivalence $s$ of $\calC$, is \emph{ample} if
\begin{enumerate}
\renewcommand{\theenumi}{B\arabic{enumi}}
\item\label{def:ample1} For all $\calM \in \calC$, there exist positive integers $l_1, \dots, l_p$
and an epimorphism $\oplus_{i=1}^p s^{-l_i}\OO \twoheadrightarrow \calM$.
\item\label{def:ample2} For all epimorphisms $\calM \twoheadrightarrow \calN$, there exists $n_0$
such that \[ \Hom(\OO, s^n \calM ) \to \Hom(\OO, s^n \calN) \] is an epimorphism
for $n \geq n_0$.
\end{enumerate}
\end{definition}

For convenience, we denote $\Hom(\OO, \calM)$ by $H^0(\calM)$.

\begin{proposition}\label{prop:ample-auto-then-noeth} {\upshape \cite[Thm.~4.5]{AZ}}
Let $(\calC, \OO, s)$ be as in Definition~\ref{def:ample-auto}. Suppose that the following conditions hold
\begin{enumerate}
\item The object $\OO$ is noetherian.
\item The ring $R_0  = H^0(\OO)$ is right noetherian and $H^0(\calM)$
is a finitely generated $R_0$-module for all $\calM \in \calC$.
\item The pair $(\OO,s)$ is ample.
\end{enumerate}
Then $R  = \Gamma(\calC, \OO, s)_{\geq 0}$ is a right noetherian
$A$-algebra and $\proj R \cong (\calC, \OO)$.\qed
\end{proposition}

Given an automorphism $\gs$ of $X$ and an invertible sheaf $\LL$, we can
define the autoequivalence $s = \LL_\gs \otimes -$ on $\coh(X)$.
For a coherent sheaf $\F$, define $\LL_\gs \otimes \F = \LL \otimes \gs^*\F$.
These are the only autoequivalences of $\coh(X)$ which we will examine,
due to the following proposition.

\begin{proposition}\label{prop:only-autos} 
{\upshape (\cite[Cor.~6.9]{AZ}, \cite[Prop.~2.15]{AV})}
Let $X$ be  proper over a field. Then any autoequivalence $s$ of $\coh(X)$ is
naturally isomorphic to $\LL_{\gs} \otimes -$ for some
automorphism $\gs$ and invertible sheaf $\LL$.\qed
\end{proposition}

Denote pull-backs by $\gs^* \F = \F^\gs$.
We then have
\[
\F \otimes (\LL_\gs)^n = \F \otimes \LL \otimes \LL^\gs \otimes \dots \otimes \LL^{\gs^{n-1}}
\]
and
\[
s^n \F = (\LL_\gs)^n \otimes \F = 
\LL \otimes \LL^\gs \otimes \dots \otimes \LL^{\gs^{n-1}} \otimes \F^{\gs^n}.
\]

\begin{definition}
Given an automorphism $\gs$ of a scheme $X$,  an invertible sheaf
$\LL$ is \emph{left $\gs$-ample} if for all coherent sheaves $\F$, there exists $n_0$
such that
\[
H^q(X, (\LL_\gs)^n \otimes  \F) = 0
\]
for $q > 0$ and $n \geq n_0$. An invertible sheaf $\LL$ is \emph{right $\gs$-ample}
if for all coherent sheaves $\F$, there exists $n_0$ such that
\[
H^q(X, \F \otimes (\LL_\gs)^n) = 0
\]
for $q > 0$ and $n \geq n_0$.
\end{definition}

\begin{lemma} {\upshape \cite[Lem.~2.4]{Ke1} }
An invertible sheaf $\LL$ is right $\gs$-ample if and only if $\LL$ is 
left $\gs^{-1}$-ample.\qed
\end{lemma}

\begin{lemma}\label{lem:chi1-then-chi}
The pair $(\OO_X, \LL_\gs \otimes -)$ is an ample autoequivalence if and 
only if $\LL$ is left $\gs$-ample.
\end{lemma}
\begin{proof}
If $s = \LL_\gs \otimes -$ is an ample autoequivalence, then the 
sequence of invertible sheaves 
\[
\{ (\gs^*)^{-n+1}(\OO_X \otimes (\LL_\gs)^n) \} \cong \{ \OO_X \otimes (\LL_{\gs^{-1}})^n \}
\] 
satisfies
condition~ \eqref{th:main2} (and hence condition~\eqref{th:main1}) 
of Theorem~\ref{th:maintheo}, so $\LL$ is right $\gs^{-1}$-ample. Thus
$\LL$ is left $\gs$-ample.

Assuming $\LL$ is left $\gs$-ample, we have that $\LL$ is right $\gs^{-1}$-ample.
So condition~\eqref{th:main1} (and hence conditions~\eqref{th:main2} and \eqref{th:main3})
of Theorem~\ref{th:maintheo} hold for the sequence $\{ \OO_X \otimes (\LL_{\gs^{-1}})^n \}$. 
This immediately gives \eqref{def:ample2}
of the definition of ample autoequivalence. 
Now because \eqref{th:main3} holds for the sequence 
$\{ \OO_X \otimes (\LL_{\gs^{-1}})^n \}$,
given a coherent sheaf $\F$,
we may pull back $\gs^*\F \otimes (\LL_{\gs^{-1}})^n$ 
by $(\gs^*)^{n-1}$ and have that for all sufficiently large $n$,
the sheaf
\[
(\LL_\gs)^n \otimes  \F = s^n \F
\]
is generated by global sections. Thus choosing one large $n_0$, we have some $p$
so that there is an epimorphism $\oplus_{i=1}^p s^{-n_0} \OO_X \twoheadrightarrow \F$.
Thus \eqref{def:ample1} of the definition of ample autoequivalence holds.
\end{proof}

The following now follows from Corollary~\ref{cor:projective}.

\begin{corollary}\label{cor:s-ample-then-projective}
Let $X$ be proper over a commutative noetherian ring $A$. Suppose that there
exists an automorphism $\gs$ and an invertible sheaf $\LL$ on $X$ such that
$\LL$ is right
 $\gs$-ample (or equivalently such
 that $\LL_\gs \otimes -$ is an ample autoequivalence). Then $X$ is projective over $A$.\qed
\end{corollary}

Let $B(X, \gs, \LL)^{\op} = \Gamma(\coh(X), \OO_X, \LL_\gs \otimes - )_{\geq 0}$. The
ring $B(X, \gs, \LL)$ is called a \emph{twisted homogeneous coordinate ring} for $X$.
Through pull-backs by powers of $\gs$, one can show the following.

\begin{lemma} {\upshape \cite[Lem.~2.2.5]{Ke-thesis} }
The rings $B(X, \gs, \LL)$ and $B(X, \gs^{-1}, \LL)$ are opposite rings.\qed
\end{lemma}

The rings $B(X, \gs, \LL)$ were extensively studied in \cite{AV}, though only
when $A$ was a field. In that paper,
the multiplication was defined using left modules instead of the right modules
used in \cite{AZ}. Thus we stipulate that $B(X, \gs, \LL)$
is the opposite ring of $\Gamma(\coh(X), \OO_X, \LL_\gs \otimes - )_{\geq 0}$,
in order to keep the multiplication consistent between \cite{AV} and
\cite{AZ}. From the previous two lemmas and
Proposition~\ref{prop:ample-auto-then-noeth} we now have the following.

\begin{proposition}
Let $X$ be proper over a commutative noetherian ring $A$. Let $\gs$ be an automorphism
of $X$ and let $\LL$ be an invertible sheaf on $X$. 
If $\LL$ is left $\gs$-ample, then $B(X, \gs, \LL)$
is left noetherian. If $\LL$ is right $\gs$-ample, 
then $B(X, \gs, \LL)$ is right noetherian.\qed
\end{proposition}

Fix an automorphism $\gs$ of $X$ and set $\LL_m = \OO_X \otimes (\LL_\gs)^m$.
For a graded ring $R = \oplus_{i=0}^\infty R_i$, the Veronese subring
$R^{(m)} = \oplus_{i=0}^\infty R_{mi}$. If $R$ is commutative and noetherian,
then some Veronese subring of $R$ is generated in degrees $0$ and $1$
\cite[p.~204, Lemma]{MuRed}; however, there are noncommutative noetherian
graded rings such that no Veronese subring is generated in
degrees $0$ and $1$ \cite[Corollary~3.2]{StafZ}.

\begin{proposition}\label{prop:veronese}
Let  $X$ be proper over a commutative noetherian ring $A$, and
let $\LL$ be $\gs$-ample on $X$. 
There exists $n$ such that
$B(X,\gs, \LL)^{(n)}$ is generated in degrees $0$ and $1$ over $A$.
\end{proposition}
\begin{proof}
Note that $B(X, \gs, \LL)^{(n)} \cong B(X, \gs^n, \LL_n)$
and $\LL$ is $\gs$-ample if and only if $\LL_n$ is $\gs^n$-ample
\cite[Lemma~4.1]{AV}.
Then by Proposition~\ref{prop:ample-poset-v-ample}, we may replace $\LL$
by  $\LL_n$ and $B$ by $B^{(n)}$ and assume $\LL$ is very ample for $\Spec A$.

For this proof, it is easiest to use the multiplication defined in \cite{AV},
namely the maps
\[
H^0(X, \LL_m) \otimes H^0(X, \LL_n) \tilde{\to} H^0(X, \LL_m) \otimes H^0(X, (\gs^*)^m\LL_n)
\to H^0(X, \LL_{m+n}).
\]
Now choose $n_0$ so that $H^q(X, \LL^{-q-1} \otimes \LL_{n}) = 0$
for $q = 1, \dots, \cd(X)$ and $n \geq n_0$. Then $\gs^*\LL_{n-1}$ is $0$-regular
with respect to $\LL$.
By Proposition~\ref{prop:regularity},
 the natural  map
\[
 H^0(X, \LL) \otimes H^0(X, \gs^* \LL_{n-1}) \to H^0(X, \LL_{n})
\]
is surjective for $n \geq n_0$.
So twisting by $(\gs^*)^i$, the natural map
\begin{equation}\label{eq:surj-map}
H^0(X, (\gs^*)^i \LL) \otimes H^0(X, (\gs^*)^{i+1}\LL_{n-1}) \to H^0(X, (\gs^*)^i \LL_{n})
\end{equation}
is surjective for $n \geq n_0$ and all $i \in \ZZ$.

Now let $\ell > 0$.
According to \eqref{eq:surj-map}, the maps
\begin{multline*} 
H^0(X, (\gs^*)^{n_0 - j}\LL) \otimes H^0(X, (\gs^*)^{n_0 - j+1} \LL_{\ell n_0 + j - 1}) \\
\to H^0(X, (\gs^*)^{n_0 - j} \LL_{\ell n_0 + j })
\end{multline*}
are surjective for $j = 1, \dots, n_0$.  
Thus the map
\[
H^0(X, \LL_{n_0}) \otimes H^0(X, (\gs^*)^{n_0} \LL_{\ell n_0} ) \to H^0(X, \LL_{(\ell + 1) n_0} ).
\]
is surjective and $B(X, \gs, \LL)^{(n_0)}$ is generated in degrees $0$ and $1$.
\end{proof}

Given an automorphism $\gs$ of $X$, let $P_\gs$ be the action of $\gs$
on $A^1(X)$. Thus $P_\gs \in \GL_\rho(\ZZ)$ for some $\rho$ by the Theorem of the Base
\ref{th:base}. We call $\gs$ quasi-unipotent if $P_\gs$ is quasi-unipotent, that is,
when all eigenvalues of $P_\gs$ are roots of unity. The statement ``$\gs$ is quasi-unipotent''
is well-defined.

\begin{proposition}
Let $X$ be proper over a noetherian ring $A$ and let $\gs$ be an
automorphism of $X$. Let $P, P' \in \GL_\rho(\ZZ)$ be two representations
of the action of $\gs$ on $A^1(X)$. Then $P$ is quasi-unipotent
if and only if $P'$ is quasi-unipotent.
\end{proposition}
\begin{proof}
Let $P$ be quasi-unipotent. We may replace $\gs$ by $\gs^i$
and assume $P = I + N$ for some nilpotent matrix $N$.
Then for all invertible sheaves $\calH$ and contracted integral curves $C$,
the function $f(m) = (\calH^{\gs^m}.C) = (P^m \calH.C)$ is a polynomial.

However, if $P'$ not quasi-unipotent, then $P'$ has an eigenvalue $r$ of absolute
value greater than $1$ \cite[Lemma~3.1]{Ke1}.
(If $X$ is projective, then the cone of nef invertible sheaves has a 
non-empty interior and
we may assume $r$ is real
 \cite[Theorem~3.1]{V}.)
Let $v = a_1 \calH_1 + \dots + a_\rho \calH_\rho \in A^1(X) \otimes \CC$ 
be an eigenvector for $r$
where the $\calH_i$ are invertible sheaves. Then
there exists a contracted integral curve $C$ such that
\[
\vert a_1 (\gs^m \calH_1 . C) \vert + \dots + \vert a_\rho (\gs^m \calH_\rho.C) \vert
\geq \vert (\gs^m v . C) \vert = \vert r \vert^m \cdot \vert (v. C) \vert > 0.
\]
Thus not all of the $(\gs^m \calH_i . C)$ can be polynomials. Thus we have
a contradiction and $P'$ must be quasi-unipotent.
\end{proof}

We now can state the following generalization of \cite[Theorem~1.3]{Ke1}.

\begin{theorem}
Let $X$ be proper over a commutative noetherian ring $A$. Let $\gs$ be an automorphism
of $X$ and let $\LL$ be an invertible sheaf on $X$. 
Then $\LL$ is right $\gs$-ample
if and only if $\gs$ is quasi-unipotent
and
\[
\LL \otimes \LL^\gs \otimes \dots \otimes  \LL^{\gs^{m-1}}
\]
is an ample invertible sheaf for some $m > 0$.
\end{theorem}
\begin{proof}
The proof mostly proceeds as in \cite[$\S$3--4]{Ke1},
using the fact that the sequence $\{ \OO_X \otimes (\LL_\gs)^m \}$ is an
ample sequence if and only if \eqref{th:main4} holds,
and then showing the equivalence of \eqref{th:main4} with the condition above. 
However, when
$\gs$ is not quasi-unipotent, one must use the method outlined in
\cite[Remark~3.5]{Ke1}, since the proof of \cite[Theorem~3.4]{Ke1}
relied on the growth of the dimensions of the graded pieces
of $B(X, \gs, \LL)$.
\end{proof}

Now since $\gs$ is quasi-unipotent if and only if $\gs^{-1}$ is quasi-unipotent,
we easily get the following, as proved in \cite[$\S$5]{Ke1}.

\begin{theorem}\label{th:thcr-noeth}
Let $A$ be a commutative noetherian ring, 
let $X$ be proper over  $A$, let $\gs$ be an automorphism
of $X$, and let $\LL$ be an invertible sheaf on $X$.
 Then $\LL$ is right $\gs$-ample if and only if $\LL$
is left $\gs$-ample. Thus we may simply say that such an $\LL$ is $\gs$-ample.
If $\LL$ is $\gs$-ample, then $B(X, \gs, \LL)$
is noetherian.\qed
\end{theorem}

\begin{remark}
Now suppose that $X$ is proper over a field $k$. Then the claims of \cite[$\S$6]{Ke1} regarding
the Gel'fand-Kirillov dimension of $B(X, \gs, \LL)$, with $\gs$-ample $\LL$,
are still valid. This is because the proofs rely on a weak Riemann-Roch formula
\cite[Example~18.3.6]{Fu} which is valid over an arbitrary field.
\end{remark}

\begin{definition} \cite[Definition~3.2]{AV}
Let $k$ be a field and
let $R$ be a finitely $\NN$-graded right noetherian $k$-algebra.
That is, $R = \oplus_{i=0}^\infty R_i$ and $\dim_k R_i$ is finite for all $i$. 
The ring $R$ is said to satisfy right $\chi_j$ if
for all finitely generated, graded right $R$-modules $M$ and all $\ell \leq j$,
\[
\dim_k \Ext^{\ell}(R/R_{> 0}, M) < \infty,
\]
where $\Ext$ is the ungraded $\Ext$-group, calculated in the category of \emph{all}
right $R$-modules.
If $R$ satisfies right $\chi_j$ for all $j \geq 0$, we say $R$ satisfies right $\chi$.
Left $\chi_j$ and left $\chi$ are defined similarly with left modules.
\end{definition}

\begin{theorem}\label{th:classify-noncomm-proj}
Let $k$ be a field and 
let $R$ be a finitely $\NN$-graded 
right noetherian $k$-algebra which satisfies right $\chi_1$.
Suppose that there exists
a scheme $X$, proper over $k$, such that $\proj R \cong (\coh(X), \OO_X)$.
Let $\rho$ be the Picard number of $X$.
Then 
\begin{enumerate}
\item\label{th:class1} $X$ is projective over $k$,
\item\label{th:class2} $R$ is noetherian and satisfies left and right $\chi$,
\item There exists $m$ such that the Veronese subring $R^{(m)}$ is generated
in degrees $0$ and $1$, and
\item $\GKdim R$ is an integer and
\[
\dim X + 1 \leq \GKdim R \leq 
 2 \left\lfloor \frac{\rho - 1}{2} \right\rfloor (\dim X - 1) + \dim X + 1.
\]
\end{enumerate}
\end{theorem}
\begin{proof}
These claims depend only on the behavior of $R$ in high degree. Thus using 
\cite[Theorem~4.5]{AZ}, we may assume $R = \Gamma(\coh(X), \OO_X, s)_{\geq 0}$ for
some autoequivalence $s$. But by Proposition~\ref{prop:only-autos}, we may assume
$s = \LL_\gs \otimes -$ for some invertible sheaf $\LL$ and automorphism $\gs$.
Thus we may assume $R^{\op} = B(X, \gs, \LL)$.

By Lemma~\ref{lem:chi1-then-chi}, the sheaf $\LL$ is $\gs$-ample 
since $s$ is an ample autoequivalence. 
So $X$ is projective by Corollary~\ref{cor:s-ample-then-projective}.
Also $R$ is noetherian by Theorem~\ref{th:thcr-noeth} and
the vanishing higher cohomology of $s^m\F$ for all coherent sheaves $\F$
gives that $R$ satisfies right $\chi$ \cite[Theorem~7.4]{AZ}.
Since $\LL$ is also $\gs^{-1}$-ample,  $R$ satisfies left $\chi$ by symmetry. 
The claim regarding
Veronese subrings is Proposition~\ref{prop:veronese}.

Finally, the claim regarding GK-dimension comes from \cite[Theorem~6.1]{Ke1}.
We need only explain the bounds. First, $\GKdim B(X,\gs,\LL) = \GKdim B(X, \gs^m, \LL_m)$
\cite[p.~263]{AV}, so we may assume that $\gs$ fixes the irreducible components of $X$.
For each irreducible component $X_i$, let $\gs_i$ be the induced automorphism.
Then \cite[Proposition~6.11]{Ke1} show that
\[
\GKdim B(X, \gs, \LL) = \max_{X_i} \GKdim B(X_i, \gs_i, \LL\vert_{X_i}).
\]
We may also assume
that $P_\gs$ is unipotent, so 
write $P_\gs = I + N$ for some nilpotent matrix $N$ and let $\ell$ be the smallest
integer so that $N^{\ell + 1} = 0$. If $P_{\gs_i} = I + N_i$, then
$N_i^{\ell + 1} = 0$, which can be seen by
pulling-back an ample invertible sheaf to $X_i$ and using \cite[Lemma~4.4]{Ke1}. 
So to find the desired bounds, we may assume $X$
is irreducible, hence equidimensional,
so the bounds in \cite[Theorem~6.1]{Ke1} apply. Now $\ell$ is even \cite[Lemma~6.12]{Ke1}
and $0 \leq \ell \leq \rho - 1$, so the universal bounds follow.
\end{proof}

Theorem~\ref{th:classify-noncomm-proj}\eqref{th:class1} seems to be a fortunate
result. It says that we cannot have a noncommutative projective scheme 
$\proj R = (\coh(X), \OO_X)$
coming from a commutative non-projective scheme $X$.

\begin{example}
Theorem~\ref{th:classify-noncomm-proj}{\eqref{th:class2}} may not be true if a structure
sheaf other than $\OO_X$ is used. In \cite[Example~4.3]{StafZ}, a coherent
sheaf $\F$ 
and ample autoequivalence $s$ 
on $\PP^1$ is chosen so that $R= \Gamma(\coh(\PP^1), \F, s)_{\geq 0}$
is right noetherian and satisfies $\chi_1$, but $R$ is not
left noetherian and does not satisfy $\chi_2$.
\end{example}

\section*{Acknowledgements}
Thanks go to M.~Artin, I.~Dolgachev, D.~Eisenbud, M.~Hochster, J.~Howald, R.~Lazarsfeld,
S.~Schr{\"o}er, J.T.~ Stafford, and B.~Totaro for useful discussions.


\providecommand{\bysame}{\leavevmode\hbox to3em{\hrulefill}\thinspace}


\newpage

\numberwithin{equation}{section}
\renewcommand{\theequation}{E\arabic{section}.\arabic{equation}}
\setcounter{section}{0}
\renewcommand*{\theHsection}{E.\the\value{section}}

\begin{center}
\LARGE Erratum
\end{center}

\section{Corrected Lemma}
\label{SecErratum}

As stated, \cite[Lemma~\ref{lem:geometrically-integral}]{E-AmpleFilters} is definitely incorrect. 
Take $A=K$ a field and $X$ a
projective scheme over $K$, with $X$ integral but not geometrically irreducible
over $K$. Let $X'$ be a scheme surjecting onto $X$. Then $X'$ cannot
be geometrically integral over $K$ as the lemma asserts as possible.
If $X'$ were geometrically integral over $K$, then $X' \times_K \overline{K}$
has only one irreducible component, while $X \times_K \overline{K}$ has
more than one. So $X' \times_K \overline{K} \to X \times_K \overline{K}$
is not surjective. But surjectivity is preserved under base change
\cite[Tag 01S1]{E-stacks-project},
yielding a contradiction.
Our mistake was treating \emph{geometrically integral} as an absolute notion,
rather than a relative one.

Before correcting the lemma, we give a presumably well-known result, due to lack
of reference. We follows the notation of \cite{E-AM}.

\begin{lemma}\label{lem:finite-opens}
Let $A, A'$ be noetherian integral domains with 
$\phi:A \to A'$ a finite ring homomorphism.
Let $S=\Spec A, S'=\Spec A'$ with induced morphism $g: S' \to S$. 
Write $S_b = S \setminus V(b)$ for $b \in A$
(and similarly for $S'_c$ with $c\in A'$). Then
for any non-zero $c \in A'$, there exists non-zero $b \in A$ such that
$S'_{\phi(b)} = g^{-1}(S_b) \subseteq S'_c $.
\end{lemma}
\begin{proof} The claim that $S'_{\phi(b)} = g^{-1}(S_b)$ is 
from \cite[Ch.~ 1, Exercise~ 21]{E-AM}. For simplicity, we can replace
$A$ with $\phi(A)$ and thus assume $A \subseteq A'$.

Choose non-zero $c \in A'$.  Since $c$ is integral over $A$, there
exists a minimal monic polynomial $f$ for $c$ with non-zero constant term $b \in A$.
This gives $b = ec$ for some $e \in A'$. Then $S'_b = S'_{ec} = S'_c \cap S'_e$
\cite[Ch. 1, Exercise~17]{E-AM}, so $S'_b \subseteq S'_c$.
\end{proof}

\begin{remark}
For the remainder of this erratum, we will abuse notation by 
replacing $S_b, S'_b$ with the isomorphic schemes $\Spec A_b, \Spec A'_b$ respectively.
\end{remark}

Here is the corrected lemma.

\begin{lemma}\label{E-lem:geometrically-integral}
Let $A$ be a noetherian domain, 
  let $X$ be an integral scheme
with  a proper, surjective morphism $\pi\colon X \to \Spec A$, 
and let $d$ be the dimension of the generic fiber of $\pi$. 
Then there exists non-zero $b \in A$, a noetherian domain $A'$,
 an integral scheme $X'$, and 
finite-type morphisms $f, g, \pi'$ with commutative diagram
\begin{equation}\label{eq:nicemorphisms}
\xymatrix{ X' \ar[r]^f \ar[d]_{\pi'} & X_b \ar[d]^{\pi_b} \\
			\Spec A' \ar[r]^g & \Spec A_b }
\end{equation}
where $X_b = X \times_A A_b$ and $\pi_b = \pi \times_A \id_{A_b}$.

The diagram has the following properties:
\begin{enumerate}
\item The morphisms $f,g,\pi'$ are projective and surjective.
\item The morphisms $g,\pi_b$ are flat.
\item\label{part3} The morphism $g$ is finite.
\item\label{part4} The morphism $\pi'$ is smooth.
\item\label{part5} Every fiber $X'_s$ of $\pi'$ is geometrically integral
over its base field $K'(s)$ and $X'_s$ has dimension $d$.
\item\label{part6} There exists an open subscheme $U \subseteq X_b$ such that
$f$ restricted to $f^{-1}(U)$ is finite.
\end{enumerate}
\end{lemma}

\begin{proof}
Let $K$ be the fraction field of $A$ and let $X_0 = X \times_A K$
be the generic fiber of $\pi$. Via Alteration of Singularities \cite[Theorem~4.1]{E-deJong}, 
there exists a finite extension $K'$ of $K$, an integral scheme $X'_0$, and commutative
diagram
\begin{equation}\label{eq:nicegeneric}
\xymatrix{ X'_0 \ar[r]^{f_0} \ar[d]_{\pi'_0} & X_0 \ar[d]^{\pi_0} \\
			\Spec K' \ar[r]^{g_0} & \Spec K }
\end{equation}
such that 
\begin{enumerate}
\item $f_0$ is an alteration (that is, a dominant, proper, generically finite
morphism), 
\item $g_0$ is finite, 
\item $g_0 \circ \pi'_0 = \pi_0 \circ f_0$ is projective,
\item $X'_0$ is geometrically irreducible and smooth over $K'$ \cite[Remark~4.2]{E-deJong}.
\end{enumerate}

Since $g_0 \circ \pi'_0 = \pi_0 \circ f_0$ is projective, we have that
$\pi'_0$ and $f_0$ are projective \cite[Exercise~II.4.8]{E-Hartshorne}. Since $f_0$ is proper
(hence closed) and dominant, we have $f_0$ surjective. Clearly the other morphisms
are surjective. Since $K$ is a field, the maps $g_0, \pi_0$ are trivially flat.
Also, since $f_0$ is generically finite, $\dim(X'_0)=\dim(X_0) = d$
\cite[2.20]{E-deJong}.

Finally, since $\pi'_0$ is finite-type, we have that $\pi'_0$ smooth
over $K'$ is equivalent to $X'_0$ geometrically regular over $K'$ 
\cite[Tag 038X]{E-stacks-project}. Thus $X'_0$ is geometrically normal
and hence geometrically reduced over $K'$ \cite[Tags 0569, 033K]{E-stacks-project}.
So $X'_0$ is geometrically integral over $K'$. 

Thus, all the properties of the lemma hold over $\Spec K$.
We need to verify that these properties can be spread out.
If $\{ A_b \}$ is the inductive system 
of one element localizations of $A$, then $K = \varinjlim A_b$
(and hence $\Spec K = \varprojlim \Spec A_b$).
By \cite[Proposition~6.2]{E-Illusie} and surrounding discussion, there exists
non-zero $b\in A$ which gives
a diagram of the form \eqref{eq:nicemorphisms} that localizes to diagram 
\eqref{eq:nicegeneric}. 

To avoid repetition, for the remainder of the proof we will state assumptions that can be
made by replacing $b$ with a multiple of $b$, and hence shrinking $\Spec A_b$,
sometimes without explicit mention of the shrinking.

\emph{A priori}, the morphisms are only of finite-type and
 the lower left and upper right corners are finite-type $A_b$-schemes, say $S',Y$
 respectively.
By the discussion in \cite{E-Illusie}, the choice of $Y$ is unique up to isomorphism
(after sufficiently shrinking $\Spec A_b$) so we may assume $Y=X \times_A A_b$
and the morphism $X \times_A A_b \to \Spec A_b$ is $\pi_b = \pi \times_A \id_{A_b}$.
Also, we may assume that $g$ is finite \cite[$\mathrm{IV}_3$, 8.10.5]{E-EGA},
so $S = \Spec A'$ where $A'$ is a finite $A_b$-module, and hence noetherian.

By Generic Freeness \cite[$\mathrm{IV}_2$, 6.9.2]{E-EGA}, we may assume that
$A'$ is a free $A_b$-module. Then $A'$ must be an integral domain; if $A'$ had a
zero-divisors $y,z$ with $yz=0$, then tensoring with $K$ would annihilate $y$ or $z$, 
contradicting
the freeness of $A'$. Similarly, $X'$ and $X_b$ have open affine
covers with coordinate rings that are free $A_b$-modules and hence integral domains.
So $X'$ is integral. This also shows that $g$ and $\pi_b$ are flat,
as this is the proof of Generic Flatness in \cite[$\mathrm{IV}_2$, 6.9.1]{E-EGA}.

We can also assume that $f, g, \pi'$ are projective and surjective 
\cite[$\mathrm{IV}_3$, 8.10.5]{E-EGA}. So we have verified the claims of the
lemma through part \eqref{part3}.

By Lemma~\ref{lem:finite-opens}, 
any open condition for $\pi'$ can be achieved by shrinking $\Spec A_b$,
as can any condition following from $\Spec K' = \varprojlim \Spec A'_b$
with $b \in A$. 
Therefore, we can assume that $\pi'$ is smooth \cite[Proposition~6.3]{E-Illusie},
that the fibers of $\pi'$ are geometrically integral over their base field
\cite[$\mathrm{IV}_3$, 12.2.4]{E-EGA}, and the fibers of $\pi'$ have
dimension $d$ \cite[$\mathrm{IV}_3$, 13.1.5]{E-EGA}. The last claim
comes from the upper semicontinuity of fiber dimension when $\pi'$ is proper.
So we have verified parts \eqref{part4} and \eqref{part5}.


Finally, let $\nu$ be the generic point of $X$. Then $\nu$ is also the generic
point of $X_0$. Since $f_0$ is generically finite, the set
$f_0^{-1}(\nu) = f^{-1}(\nu)$ is finite. Then $f$ is generically finite
and \eqref{part6} holds
by \cite[Exercise~II.3.7]{E-Hartshorne}.
\end{proof}

\begin{remark}
Note that the corrected lemma holds regardless if the quotient field $K$ of $A$
is perfect.
However, if $K$ is perfect, then $g_0 \circ \pi'_0$ is also smooth
\cite[Remark~4.2]{E-deJong}. In this case, we can also assume $g \circ \pi'$ is smooth.
\end{remark}

\begin{remark}
The original \cite[Lemma~\ref{lem:geometrically-integral}]{E-AmpleFilters} did not have a claim similar to
\eqref{part6}. While not necessary for the rest of the paper, it relativizes 
the concept of \emph{alteration} and perhaps will have a future use.

The original lemma also assumed $\pi$ projective, but assuming $\pi$ proper
required no extra work. Further, it was only argued that the generic fiber
of $\pi'$ has dimension $d$.
\end{remark}

\section{Verifying dependent results}

In this section, we verify that theorems relying on \cite[Lemma~\ref{lem:geometrically-integral}]{E-AmpleFilters}
are correct. 
There are three such theorems, 
namely \cite[Theorems~\ref{th:general-serre1}, \ref{th:base}, Proposition~\ref{prop:general-serre-integer}]{E-AmpleFilters}. 
Fortunately, having $g$ finite and faithfully flat (that is, flat and surjective) 
is close enough
to $g$ being the identity morphism.

We start with a slight generalization of \cite[Lemma~\ref{lem:surj-curve}]{E-AmpleFilters}.
Recall that for a proper $\pi: X \to S$, a closed subscheme $V \subseteq X$
is \emph{$\pi$-contracted} if $\pi(V)$ is $0$-dimensional
\cite[Definition~\ref{def:pi-contracted}]{E-AmpleFilters}.

\begin{lemma}\label{lem:pi-contracted}
Let $S, S'$ be noetherian schemes. Consider the commutative diagram
\[
\xymatrix{ X' \ar[r]^f \ar[d]_{\pi'} & X \ar[d]^\pi \\
			S' \ar[r]^g & S }
\]
with  proper morphisms $\pi, \pi'$, (proper) surjective $f$, and finite $g$. 
Let $C \subseteq X$ be a $\pi$-contracted integral curve. Then there exists a
$\pi'$-contracted integral
curve $C' \subseteq X'$ such that $f(C') = C$.
\end{lemma}
\begin{proof} By \cite[Lemma~\ref{lem:surj-curve}]{E-AmpleFilters}, there exists a $g\circ \pi'$-contracted
integral curve $C'$ such that $f(C')=C$. Since $g$ is finite, the set
$g^{-1}((g\circ \pi')(C'))$ is finite. Thus $C'$ must be $\pi'$-contracted,
as desired.
\end{proof}

We immediately have

\begin{theorem}\label{thm:36works}
Consider a diagram of schemes as in \eqref{eq:nicemorphisms}, 
satisfying the conclusions of Lemma~\ref{E-lem:geometrically-integral}.
Then the conclusions of \cite[Lemmas~\ref{lem:nef-change-scheme}, \ref{lem:nef-inject}]{E-AmpleFilters} hold.
Therefore \cite[Theorem~\ref{th:base}]{E-AmpleFilters} is true.
\end{theorem}
\begin{proof}
By Lemma~\ref{lem:pi-contracted}, the results  \cite[Lemmas~\ref{lem:nef-change-scheme}, \ref{lem:nef-inject}]{E-AmpleFilters}
hold when $g$ is finite.

The original proof of \cite[Theorem~\ref{th:base}]{E-AmpleFilters} uses
\cite[Lemma~\ref{lem:nef-inject}]{E-AmpleFilters} to replace $\pi:X_b \to \Spec A_b$ with 
$g \circ \pi': X' \to \Spec A_b$. Since \cite[Lemma~\ref{lem:nef-inject}]{E-AmpleFilters}
holds in the context of Lemma~\ref{E-lem:geometrically-integral},
we can actually replace $\pi:X_b \to \Spec A_b$
with $\pi': X' \to \Spec A'$. The remainder of the
proof of \cite[Theorem~\ref{th:base}]{E-AmpleFilters} proceeds as before,
since we can assume $\pi'$ is smooth (hence flat) and has geometrically
integral fibers.
\end{proof}

We will need \cite[Lemma~\ref{lem:nef-change-scheme}]{E-AmpleFilters} below.

The use of Lemma~\ref{E-lem:geometrically-integral} to prove 
\cite[Theorem~ \ref{th:general-serre1}, Proposition~\ref{prop:general-serre-integer}]{E-AmpleFilters}
requires examination of the technical \cite[Notation~\ref{notation:Vq}]{E-AmpleFilters}.
We repeat it here, with an extra indication of the base ring.

\begin{notation}\label{E-notation:Vq}
Let $A$ be a noetherian domain, let $\pi\colon X \to \Spec A$
be a projective morphism, let
$\LL$ be an ample invertible sheaf on $X$,  and let $\Lambda$ be a set
of (isomorphism classes of) invertible sheaves on $X$.
(For a locally closed subscheme $Y \subseteq X$, 
let $\Lambda\vert_Y = \{ \calN\vert_Y \colon \calN \in \Lambda \}$.)
Let $\F$ be a coherent sheaf on $X$.
For $q > 0$, 
we say that $V^q_A(\F, \LL, \Lambda)$ holds if there exists $m_0 = m(\F,q)$
and a non-empty open subscheme $U = U(\F, q) \subseteq \Spec A$ such that
$R^{s}\pi_*( \F \otimes \LL^m \otimes \calN)\vert_V = 0$ for all
$m \geq m_0$, $s \geq q$, $\calN \in \Lambda$, and all open subschemes $V \subseteq U$ such that 
$\calN\vert_{\pi^{-1}(V)}$ is $\pi\vert_{\pi^{-1}(V)}$-nef. 
If $\LL$ or $\Lambda$  are clear, we write $V^q_A(\F)$.
\end{notation}

We need to show that the definition of $V^q_A(\F)$ does not, in a sense, depend on $A$.

\begin{lemma}\label{lem:Vq-differentA}
Let $A,A'$ be noetherian domains, let $\pi:X \to \Spec A'$ be projective and surjective,
and  let $g: \Spec A' \to \Spec A$ be finite and surjective.
Let
$\LL$ be an ample invertible sheaf, let $\Lambda$ be a set
of invertible sheaves, and let $\F$ be a coherent sheaf on $X$. Then for $q > 0$,
if $V^q_{A'}(\F, \LL, \Lambda)$, then $V^q_A(\F, \LL, \Lambda)$.
\end{lemma}
\begin{proof}
Suppose $V^q_{A'}(\F, \LL, \Lambda)$ holds. Then there exists $m_0$
and
a non-empty open subscheme $U' \subseteq \Spec A'$
such that
$R^{s}\pi_*( \F \otimes \LL^m \otimes \calN)\vert_{V'} = 0$ for all
$m \geq m_0$, $s \geq q$, $\calN \in \Lambda$, and all open subschemes $V' \subseteq U'$ such that 
$\calN\vert_{\pi^{-1}(V')}$ is $\pi\vert_{\pi^{-1}(V')}$-nef.

Since $g$ is surjective and $A$ is reduced, 
the associated ring homomorphism $A \to A'$ is injective
\cite[Ch. 1, Exercise 21]{E-AM}, so we can assume $A \subseteq A'$. 
By Lemma~\ref{lem:finite-opens}, 
we may choose non-zero $c \in A$ with $\Spec A'_c \subseteq U'$. 
Recall that projective, surjective, and finite morphisms are preserved under base change
\cite[Tag 02V6, 01S1, 01WL]{E-stacks-project}, 
as is relative nefness \cite[Lemma~\ref{lem:field-extension}]{E-AmpleFilters}. The
vanishing condition is also local.
So we can replace $\Spec A, \Spec A',X$ with 
$\Spec A_c, \Spec A'_c, f^{-1}(\Spec A'_c)$ respectively.
Thus we may assume that $U' = \Spec A'$.

Let $m \geq m_0, s \geq q$, and $\calN \in \Lambda$. Suppose $V \subseteq \Spec A$
such that $\calN\vert_{(g\circ \pi)^{-1}(V)}$ is 
$(g \circ \pi)\vert_{(g\circ \pi)^{-1}(V)}$-nef. Then $\calN\vert_{(g\circ \pi)^{-1}(W)}$ 
is also relatively
nef, where $W$ is any affine open subscheme of $V$ \cite[Corollary~\ref{cor:nef-open}]{E-AmpleFilters}. 
But over an affine base, nefness is an absolute notion 
\cite[Proposition~\ref{prop:affine-base-not-relative}]{E-AmpleFilters}. That is, we can simply say that
$\calN\vert_{(g\circ \pi)^{-1}(W)} = \calN\vert_{\pi^{-1}(g^{-1}(W))}$ is nef.
So this invertible sheaf is also $\pi\vert_{\pi^{-1}(g^{-1}(W))}$-nef.

Set
\[
\calM' = R^{s}\pi_*( \F \otimes \LL^m \otimes \calN)\vert_{g^{-1}(W)},
\qquad \calM = R^{s}(g \circ \pi)_*( \F \otimes \LL^m \otimes \calN)\vert_{W}.
\]
Thus the $V^q_{A'}(\F, \LL, \Lambda)$ property gives 
$\calM' = 0$.
Since $g$ is finite, the scheme $g^{-1}(W)$ is affine.
Thus $\calM'$ equals the sheafification over $g^{-1}(W)$ of 
the abelian group $M$ \cite[Proposition~III.8.5]{E-Hartshorne} where
\[
M=H^s( \pi^{-1}(g^{-1}(W)), (\F \otimes \LL^m \otimes \calN)\vert_{\pi^{-1}(g^{-1}(W))}).
\]
So $M=0$.
However, $W$ is also affine, so $\calM$ is also the sheafification of $M=0$.
Thus $\calM = 0$. 

Since $\calM = 0$ holds for any affine $W \subseteq V$, we have 
$R^{s}(g \circ \pi)_*( \F \otimes \LL^m \otimes \calN)\vert_{V} = 0$ by
definition of right derived functor. Hence $V^q_A(\F, \LL, \Lambda)$
holds.
\end{proof}

\begin{theorem}
Lemma~\ref{E-lem:geometrically-integral} can be used to complete the reductions
which help prove \cite[Theorem~\ref{th:general-serre1}, Proposition~\ref{prop:general-serre-integer}]{E-AmpleFilters}.
\end{theorem}
\begin{proof}
Let $A$ be a noetherian ring and $\pi: X \to \Spec A$ be a projective morphism.
Let $\LL$ be an ample invertible sheaf and $\Lambda$ a certain set of 
invertible sheaves on $X$.
The goal of both of \cite[Theorem~\ref{th:general-serre1}, Proposition~\ref{prop:general-serre-integer}]{E-AmpleFilters}
is to show $V^1_A(\F, \LL, \Lambda)$ holds for all coherent sheaves $\F$.
Standard reductions and the definition of $V^1_A$ allows one
to assume that $\pi$ is surjective and that $X$ and $\Spec A$ are
integral. (See the original proofs for details.)

Then \cite[Lemma~\ref{lem:cover-reduction2}]{E-AmpleFilters}
and the definition of $V^1_A$
allow the reduction to $g \circ \pi':X' \to \Spec A_b$ for some $b \in A$,
where $\pi'$ has the nice properties of Lemma~\ref{E-lem:geometrically-integral}.

Finally, Lemma~\ref{lem:Vq-differentA} allows the reduction to
$\pi': X' \to \Spec A'$. The original proofs then proceed as before.
\end{proof}

\section*{Acknowledgement}
We thank R.~Lazarsfeld for pointing out the error.






\end{document}